\documentclass[11pt]{article}

\usepackage{amssymb}

\newcommand{\newc}{\newcommand}

\newc{\eqnoset}{\setcounter{equation}{0}}

\newcommand{\mref}[1]{(\ref{#1})}
\newcommand{\reflemm}[1]{Lemma~\ref{#1}}
\newcommand{\refrem}[1]{Remark~\ref{#1}}
\newcommand{\reftheo}[1]{Theorem~\ref{#1}}

\newcommand{\refcoro}[1]{Corollary~\ref{#1}}
\newcommand{\refprop}[1]{Proposition~\ref{#1}}
\newcommand{\refsec}[1]{Section~\ref{#1}}

\newcommand{\beq}{\begin{equation}}
\newcommand{\eeq}{\end{equation}}
\newcommand{\beqno}[1]{\begin{equation}\label{#1}}

\newcommand{\barr}{\begin{array}}
\newcommand{\earr}{\end{array}}

\newc{\bearr}{\begin{eqnarray*}}
\newc{\eearr}{\end{eqnarray*}}

\newc{\bearrno}[1]{\begin{eqnarray}\label{#1}}
\newc{\eearrno}{\end{eqnarray}}

\newc{\non}{\nonumber}
\newc{\nol}{\nonumber\nl}

\newcommand{\bdes}{\begin{description}}
\newcommand{\edes}{\end{description}}
\newc{\benu}{\begin{enumerate}}
\newc{\eenu}{\end{enumerate}}
\newc{\btab}{\begin{tabular}}
\newc{\etab}{\end{tabular}}

\newtheorem{theorem}{Theorem}[section]
\newtheorem{defi}[theorem]{Definition}
\newtheorem{lemma}[theorem]{Lemma}
\newtheorem{rem}[theorem]{Remark}
\newtheorem{exam}[theorem]{Example}
\newtheorem{propo}[theorem]{Proposition}
\newtheorem{corol}[theorem]{Corollary}

\newcommand{\btheo}[1]{\begin{theorem}\label{#1}}
\newc{\brem}[1]{\begin{rem}\label{#1}\em}
\newc{\bexam}[1]{\begin{exam}\label{#1}\em}
\newc{\bdefi}[1]{\begin{defi}\label{#1}}
\newcommand{\blemm}[1]{\begin{lemma}\label{#1}}
\newcommand{\bprop}[1]{\begin{propo}\label{#1}}
\newcommand{\bcoro}[1]{\begin{corol}\label{#1}}
\newcommand{\etheo}{\end{theorem}}
\newcommand{\elemm}{\end{lemma}}
\newcommand{\eprop}{\end{propo}}
\newcommand{\ecoro}{\end{corol}}
\newc{\erem}{\end{rem}}
\newc{\eexam}{\end{exam}}
\newc{\edefi}{\end{defi}}

\newc{\rmk}[1]{{\bf REMARK #1: }}
\newc{\DN}[1]{{\bf DEFINITION #1: }}

\newcommand{\bproof}{{\bf Proof:~~}}
\newc{\eproof}{{\vrule height8pt width5pt depth0pt}\vspace{3mm}}

\newc{\bfrac}[2]{\dspl{\frac{#1}{#2}}}

\newc{\nid}{\noindent}

\newcommand{\dspl}{\displaystyle}
\newc{\grad}{\nabla}
\newc{\Div}{\mbox{div}}
\newc{\pdt}[1]{\dspl{\frac{\partial{#1}}{\partial t}}}
\newc{\pdn}[1]{\dspl{\frac{\partial{#1}}{\partial \nu}}}
\newc{\pdNi}[1]{\dspl{\frac{\partial{#1}}{\partial \mathcal{N}_i}}}
\newc{\pD}[2]{\dspl{\frac{\partial{#1}}{\partial #2}}}
\newc{\dt}{\dspl{\frac{d}{dt}}}
\newc{\bdry}[1]{\mbox{$\partial #1$}}
\newc{\sgn}{\mbox{sign}}

\newc{\Hess}[1]{\frac{\partial^2 #1}{\pdh z_i \pdh z_j}}
\newc{\hess}[1]{\partial^2 #1/\pdh z_i \pdh z_j}

\newc{\ag}{\alpha}
\newc{\bg}{\beta}
\newc{\cg}{\gamma}\newc{\Cg}{\Gamma}
\newc{\dg}{\delta}\newc{\Dg}{\Delta}
\newc{\eg}{\varepsilon}
\newc{\zg}{\zeta}
\newc{\thg}{\theta}
\newc{\llg}{\lambda}\newc{\LLg}{\Lambda}
\newc{\kg}{\kappa}
\newc{\rg}{\rho}
\newc{\sg}{\sigma}\newc{\Sg}{\Sigma}
\newc{\tg}{\tau}
\newc{\fg}{\phi}\newc{\Fg}{\Phi}
\newc{\vfg}{\varphi}
\newc{\og}{\omega}\newc{\Og}{\Omega}
\newc{\pdh}{\partial}

\newc{\ccG}{{\cal G}}

\newc{\ii}[1]{\int_{#1}}
\newc{\iidx}[2]{{\dspl\int_{#1}~#2~dx}}
\newc{\bii}[1]{{\dspl \ii{#1} }}
\newc{\biii}[2]{{\dspl \iii{#1}{#2} }}
\newc{\su}[2]{\sum_{#1}^{#2}}
\newc{\bsu}[2]{{\dspl \su{#1}{#2} }}

\newc{\biiom}[1]{{\dspl\int_{\bdrom}~ #1 ~d\sg}}
\newc{\io}[1]{{\dspl\int_{\Og}~ #1 ~dx}}
\newc{\bio}[1]{{\dspl\int_{\bdrom}~ #1 ~d\sg}}
\newc{\bsir}{\bsu{i=1}{r}}
\newc{\bsim}{\bsu{i=1}{m}}

\newc{\iibr}[2]{\iidx{\bprw{#1}}{#2}}
\newc{\Intbr}[1]{\iibr{R}{#1}}
\newc{\intbr}[1]{\iibr{\rg}{#1}}
\newc{\intt}[3]{\int_{#1}^{#2}\int_\Og~#3~dxdt}

\newc{\itQ}[2]{\dspl{\int\hspace{-2.5mm}\int_{#1}~#2~dz}}
\newc{\mitQ}[2]{\dspl{\rule[1mm]{4mm}{.3mm}\hspace{-5.3mm}\int\hspace{-2.5mm}\int_{#1}~#2~dz}}
\newc{\mitQQ}[3]{\dspl{\rule[1mm]{4mm}{.3mm}\hspace{-5.3mm}\int\hspace{-2.5mm}\int_{#1}~#2~#3}}

\newc{\mitx}[2]{\dspl{\rule[1mm]{3mm}{.3mm}\hspace{-4mm}\int_{#1}~#2~dx}}
\newc{\mitmu}[2]{\dspl{\rule[1mm]{3mm}{.3mm}\hspace{-4mm}\int_{#1}~#2~d\mu}}
\newc{\iidmu}[2]{{\dspl\int_{#1}~#2~d\mu}}

\newc{\iidm}[3]{{\dspl\int_{#1}~#2~d #3}}

\newc{\itQmu}[2]{\dspl{\int\hspace{-2.5mm}\int_{#1}~#2~d\mu}}
\newc{\mitQmu}[2]{\dspl{\rule[1mm]{4mm}{.3mm}\hspace{-5.3mm}\int\hspace{-2.5mm}\int_{#1}~#2~d\mu}}

\newc{\mitQq}[2]{\dspl{\rule[1mm]{4mm}{.3mm}\hspace{-5.3mm}\int\hspace{-2.5mm}\int_{#1}~#2~d\bar{z}}}
\newc{\itQq}[2]{\dspl{\int\hspace{-2.5mm}\int_{#1}~#2~d\bar{z}}}

\newc{\pder}[2]{\dspl{\frac{\partial #1}{\partial #2}}}

\newc{\bdrom}{\bdry{\Og}}

\newc{\bilhom}{\mbox{Bil}(\mbox{Hom}(\RR^{nm},\RR^{nm}))}
\newc{\VV}[1]{{V(Q_{#1})}}

\newc{\ccA}{{\mathcal A}}
\newc{\ccB}{{\mathcal B}}
\newc{\ccC}{{\mathcal C}}
\newc{\ccD}{{\mathcal D}}
\newc{\ccE}{{\mathcal E}}
\newc{\ccH}{\mathcal{H}}
\newc{\ccF}{\mathcal{F}}
\newc{\ccI}{{\mathcal I}}
\newc{\ccJ}{{\mathcal J}}
\newc{\ccK}{{\mathcal K}}
\newc{\ccP}{{\mathcal P}}
\newc{\ccQ}{{\mathcal Q}}
\newc{\ccR}{{\mathcal R}}
\newc{\ccS}{{\mathcal S}}
\newc{\ccT}{{\mathcal T}}
\newc{\ccX}{{\mathcal X}}
\newc{\ccY}{{\mathcal Y}}
\newc{\ccZ}{{\mathcal Z}}

\newc{\bb}[1]{{\mathbf #1}}

\newc{\myprod}[1]{\langle #1 \rangle}
\newc{\mypar}[1]{\left( #1 \right)}

\newc{\BLLg}{\mathbf{\LLg}}

\newc{\mA}{\mathbf{A}}
\newc{\mB}{\mathbf{B}}
\newc{\mC}{\mathbf{C}}
\newc{\mD}{\mathbf{D}}
\newc{\mE}{\mathbf{E}}
\newc{\mF}{\mathbf{F}}
\newc{\mJ}{\mathbf{J}}
\newc{\mG}{\mathbf{G}}
\newc{\mP}{\mathbf{P}}
\newc{\mR}{\mathbf{R}}
\newc{\mQ}{\mathbf{Q}}
\newc{\mX}{\mathbf{X}}
\newc{\muu}{\mathbf{u}}
\newc{\mvv}{\mathbf{v}}

\newc{\mllg}{\mathbb{\lambda}}
\newc{\mLLg}{\mathbf{\LLg}}

\newc{\lspn}[2]{\mbox{$\| #1\|_{\Lsp{#2}}$}}
\newc{\Lpn}[2]{\mbox{$\| #1\|_{#2}$}}
\newc{\Hn}[1]{\mbox{$\| #1\|_{H^1(\Og)}$}}

\newc{\mynorm}[2]{\| #1\|_{#2}}

\newcommand{\RR}{{\rm I\kern -1.6pt{\rm R}}}

\newc{\itQQ}[2]{\dspl{\int_{#1}#2\,dz}}
\newc{\mmitQQ}[2]{\dspl{\rule[1mm]{4mm}{.3mm}\hspace{-4.3mm}\int_{#1}~#2~dz}}
\newc{\MmitQQ}[2]{\dspl{\rule[1mm]{4mm}{.3mm}\hspace{-4.3mm}\int_{#1}~#2~d\mu}}

\newc{\MUmitQQ}[3]{\dspl{\rule[1mm]{4mm}{.3mm}\hspace{-4.3mm}\int_{#1}~#2~d#3}}
\newc{\MUitQQ}[3]{\dspl{\int_{#1}~#2~d#3}}

\newc{\mccP}{\mathbb{P}}
\newc{\mccK}{\mathbb{K}}

\newc{\DKTmU}{\mccK(U)}
\newc{\DKTmUold}{(K_U(U)^{-1})^T}

\newc{\myPi}{\mathbf{W}}
\newc{\myIbar}{\bar{\ccI}_1}
\newc{\myIhat}{\hat{\ccI}_1}
\newc{\myIbreve}{\breve{\ccI}_0}

\newc{\mmk}{\mathbf{k}}

\newc{\mfu}{\mathbf{f_u}}
\newc{\mh}{\mathbf{h}}
\newcommand{\barrl}[2]{\barr{ll}\lefteqn{#1}\hspace{#2}&\\}

\begin{document}

\vspace*{-.8in}
\begin{center} {\LARGE\em Existence and Uniqueness of Weak Solutions to a Class of Degenerate Cross Diffusion Systems.}

 \end{center}

\vspace{.1in}

\begin{center}

{\sc Dung Le}{\footnote {Department of Mathematics, University of
Texas at San
Antonio, One UTSA Circle, San Antonio, TX 78249. {\tt Email: Dung.Le@utsa.edu}
}}

\end{center}

\begin{abstract}
We consider a class of cross diffusion systems with degenerate (or porous media type) diffusion which is  inspired by models in mathematical biology/ecology with zero self diffusions. Known techniques for scalar equations are no longer available here as maximum/comparison principles are generally unavailable for systems. However, we will provide the existence of  weak solutions to the degenerate systems under mild integrability conditions of strong solutions to nondegenerate systems and show that they converge to a weak solution of the degerate system. These conditions will be verified for the model introduced by Shigesada {\it et al.} in \cite{SKT}. Uniqueness of limiting and unbounded weak solutions will also be proved.    \end{abstract}

\vspace{.2in}

\section{Introduction}\label{introsec}\eqnoset

In this paper, we study the solvability of  the following parabolic system of $m$ equations ($m\ge2$)  
\beqno{ep1}u_t-\Delta(P(u))=f(u),\quad (x,t)\in \Og\times(0,T_0)\eeq  for the unknown vector $u=[u_i]_{i=1}^m$.
Here,  $P$ and $f$ are $C^1$ maps on $\RR^m$. 

The system is equipped with boundary and initial conditions
$$\left\{\barr{l} \mbox{$u=0$  on $\partial \Og\times(0,T_0)$},\\ u(x,0)=u_0(x) ,\quad x\in \Og. \earr\right.$$

The consideration of \mref{ep1} is motivated by the extensively studied porous media equation for a scalar unknown $u:\Og\times(0,T_0)\to\RR$ and some $k>0$ $$u_t-\Delta(|u|^ku)=f(u),\quad \mbox{ $(x,t)\in \Og\times(0,T_0)$ }.$$ There is a vast literature on this equation but, to the best of our knowledge, no work has discussed its vectorial cases. Naturally, the vectorial version of this equation is the system \mref{ep1} of $m$ equations with $P(u)=|u|^ku$,  where $|u|=\sqrt{u_1^2+\cdots+u_m^2}$.

Let $A(u)=P_u(u)$, the Jacobian of $P$, and $\llg(u)=|u|^k$. We easily see that  $$\llg(u)\le\myprod{A(u)\zeta,\zeta} \mbox{ and } |A(u)|\le (1+k)\llg(u)\quad \forall u\in\RR^m, \zeta\in\RR^{mN}.$$
This naturally leads us to the consideration the following general and main condition for the system \mref{ep1}.

\bdes

\item[P)] $P:\RR^m\to\RR^m$ is a $C^1$ map.  The Jacobian $A(u)=P_u(u)$ satisfies: there are a constant $C_*>0$  and a {\em nonnegative} scalar $C^1$ function $\llg(u)$ on $\RR^m$  such that  for all $u\in\RR^m$,  $\zeta\in\RR^{mn}$ 
\beqno{A1} \llg(u)|\zeta|^2 \le \myprod{A(u)\zeta,\zeta} \mbox{ and } |A(u)|\le C_*\llg(u).\eeq

In addition, $\llg (u)$ has a polynomial growth in $|u|$. That is,  $\llg (u)\sim|u|^k$ for some $k>0$.
\edes

Of course, the polynomial growth of $A(u), \llg (u)$ implies that $|A_u(u)|\le C|\llg_u(u)|$.

Under this assumption, \mref{ep1} is a {\em strongly coupled} parabolic system, as the matrix $A(u)$ is a full matrix in general. Moreover, we assume only that $\llg(u)\ge0$ so that $A(u)$ can be degenerate, i.e., $\llg (u)\equiv0$, in a set of $\RR^m$, say $\{0\}$. Therefore, when we discuss the existence of strong solutions to the nondegenerate system \mref{ep1} we also need to consider the following hypothesis. 
\bdes\item[PR)] P) holds and there is some $\llg_0>0$ such that $\llg(u)\ge\llg_0$ for all $u\in \RR^m$.
\edes 

For  $\llg (u)$ with polynomial growth this condition is equivalent to the assumption that $\llg (u)\sim (\llg_0+|u|)^k$ for some $k,\llg_0>0$.

Concerning the reaction term $f$,  we assume the following condition.
\bdes \item[f)] $f:\RR^m\to\RR^m$  is a $C^1$ map and there exists a constant $C$ such that 
\beqno{fumore}|f(u)|\le C|u|(1+\llg(u)),\eeq
\beqno{fxugrowth} |f_u(u)|\le C(1+\llg(u)).\eeq
\edes

The structural conditions P), in particular PR), and f) are also motivated by the well known SKT model introduced by Shigesada {\it et al.} in \cite{SKT}
\beqno{e0}\left\{\barr{lll} (u_1)_t &=& \Delta(d_1u_1+\ag_{11}u_1^2+\ag_{12}u_1u_2)+\Div[ b_1u_1\nabla\Fg(x)]+f_1(u_1,u_2),\\(u_2)_t &=& \Delta(d_2u_2+\ag_{21}u_1u_2+\ag_{22}u_2^2)+\Div[ b_2u_2\nabla\Fg(x)]+f_2(u_1,u_2).\earr\right.\eeq Here, $f_i(u_1,u_2)$ are reaction terms of Lotka-Volterra type and quadratic in $u_1,u_2$.  Dirichlet or Neumann boundary conditions were usually assumed for \mref{e0}. This model was used to describe the population dynamics of the species densities  $u,v$ which move under the influence of population pressures and the environmental potential $\Fg(x)$.

If $d_1,d_2$ are positive, under the following assumption on the constant parameters $\ag_{ij}$'s \beqno{agcondyagi} \ag_{ij}>0,\; \ag_{21}^2<8\ag_{11}\ag_{12},\; \ag_{12}^2<8\ag_{22}\ag_{21},\eeq
and that $\Og$ is a planar domain ($N=2$), Yagi proved in \cite{yag} the global existence of positive solutions, with positive initial data. In this paper, we will extend this result and related others by considering a much more general structural conditions like PR) and f). Indeed, we will replace the quadratics in the Laplacians and $f_i$ of \mref{e0} by appropriate polynomials of order $k+1$ for some $k>0$. Obviously,
the SKT system \mref{e0} is a special case of \mref{ep1} with
$P:\RR^2\to\RR^2$ being a quadratic map which
satisfies PR) for $\llg (u)$ being some linear function in $|u|$, $u=[u_1,u_2]^T$. Because $f_i$'s in \mref{e0} are quadratic in $u_1,u_2$, it is clear that  the condition f) is also verified here.

Here, we will discuss the existence of weak solutions to \mref{e0} when the self diffusion coefficients $d_1,d_2$ are zero. This is just a special case of the condition P) considered here. Again, our work may be the first addressing such problem in this general setting.

In particular, a simple consequence of our main results applying to the degenerate \mref{e0} ($d_1=d_2=0$) with Lotka-Volterra type reaction terms on planar domains ($N=2$)
\beqno{iSKT}\left\{\barr{lll}(u_1)_t&=&\Delta(u_1[\ag_{11}u_1+\ag_{12}u_2])+u_1(a_1+b_1u_1+c_1u_2),\\(u_2)_t&=&\Delta(u_2[\ag_{21}u_1+\ag_{22}u_2])+u_2(a_2+b_2u_1+c_2u_2).\earr\right.\eeq

Of course, this system is a special case of \mref{ep1} with $u=[u_1,u_2]^T$ and
$$P(u)=[u_1(\ag_{11}u_1+\ag_{12}u_2),u_2(\ag_{21}u_1+\ag_{22}u_2)]^T,$$
 $$f(u)=[u_1(a_1+b_1u_1+c_1u_2),u_2(a_2+b_2u_1+c_2u_2)]^T.$$

In literature, the system \mref{iSKT} is said to be competitive if the constants $b_i,c_i$  are nonpositive. Clearly,
\beqno{SKTcompeteintro} \myprod{f(u),u} \le C|u|^2,\;\forall u=[u_1,u_2]^T,\;u_1,u_2\ge0.\eeq

In general, we assume that there are $C_0,c_0>0$ such that \beqno{SKTcoopintro} \myprod{f(u),u} \le C_0|u|^2+c_0|u|^3,\;\forall u\in \RR^2.\eeq

We have the following easy (and new) consequence of our main results.

\bcoro{degSKT} Assume $N=2$, \mref{agcondyagi} and nonnegative initial data. Suppose further that either \bdes \item[a)]\mref{SKTcompeteintro} holds (i.e., the system is competitive);\edes or \bdes\item[b)]\mref{SKTcoopintro} holds and either that homogeneous Dirichlet boundary condition is assumed and $c_0$ is small or that homogeneous Neumann boundary condition is assumed and $\|u\|_{L^1(Q)}$ is uniformly bounded for any strong solution $u$ to \mref{e0}.\edes   Then there is a nonnegative weak solution $u=[u_1,u_2]^T$ to the degenerate system \mref{iSKT}. This solution is  the limit of strong solutions to the nondegenerate systems \mref{e0} when $d_1,d_2$ tend to $0$. Moreover, this weak solution is VMO. \ecoro

The proof of this result will be presented in \refsec{planar}.

We organize our paper as follows. In \refsec{compsec}, we first collect some basic compactness results and basic inequalities which will be used throughout this paper. We will discuss in \refsec{res} the existence of {\em strong solutions} of \mref{ep1} when it is regular (i.e., PR) holds). These results are just simple consequences of the theory for general strongly coupled parabolic and elliptic systems in \cite{dlebook}, which provides an alternative approach to the existence of strong solutions in \cite{Am2}. The results in \refsec{res} hold under very weak integrability assumptions and the crucial (but weakest) condition that the strong solutions have apriori small BMO (Bounded Mean Oscillation)  norms in small balls. Again, we would like to emphasize that no boundedness of solutions will be assumed here because maximum or comparison principles are not available for systems.

Once the existence of strong solutions for regular systems is proved, we will follow the standard approach to establish the existence of {\em weak solutions} to the degenerate systems (i.e. P) holds but $\llg (u)$ can be zero on some  subset of $\RR^m$). We will approximate the degenerate systems by a sequence of regular ones whose strong solutions can be estimated  uniformly so that we can pass to the limit, using a compactness result in \refsec{compsec}. Uniform estimates of strong solutions to \mref{ep1} will be the most important matter of this paper and they will be established in \refsec{weaksol} under very mild uniform integrability assumptions on the strong solutions and data of the approximation systems. Examples, for planar domains when these assumption can be verified, are presented in \refsec{planar} where we also provide the proof of \refcoro{degSKT}. 

In \refsec{uniweaksec} we prove the uniqueness of the limit weak solution obtained in \refsec{weaksol}. For degenerate scalar equations this has been done in several works, starting with the work of Br\'ezis and Crandall \cite{BC} which  relies on maximum/comparison principles which are not available here for systems (see also the excellent monograph \cite{Vas} on this matter).  We end the paper by establishing in \refsec{uniboundw} the uniqueness of {\em unbounded} weak solutions to {\em nondegenerate} cross diffusion systems. The class of weak solutions we consider here is much broader than those usually used in literature.

\section{Some technical lemmas} \eqnoset \label{compsec}

 We first have the following compactness result which is an improved version of \cite[Lemma 3.3]{letrans} and more suitable for our purposes here.  In the sequel and throughout this paper, we will denote by $v_t,Dv$ the temporal and spatial partial derivatives of a function $v$.

\blemm{compactlemm} Let $Q=B\times[-1,0]$. Consider sequences of functions $\{v_k\}$ on $Q$ and assume that 
\bdes \item[c.1)] There is a constant $M$  such that for all $k$ $$ \|v_k\|_{L^1(Q)},\,\|Dv_k\|_{L^2(Q)} \le M.$$
\item[c.2)] For any given $\mu>0$ there is $C(\mu)$ such that if $-1<s<r<0$ and $r-s<C(\mu)$ then $$\int_s^r\iidx{B}{|(v_k)_t|}\, d\tau\le \mu \quad \forall k.$$

\edes

Then for any  $q\in[1,\infty)$ and  $p\in [1,2_*)$ (as usual, $2_*$ is the Sobolev conjugate of $2$, i.e., it can be any number in $[1,\infty)$ if $N=2$ and $2_*=\frac{2N}{N-2}$ otherwise) the sequence $\{v_k\}$ is precompact in $L^q((-1,0),L^p(B))$.
\elemm

\bproof  First of all, if we use the equivalent norm $\|v\|_{W^{1,2}(B)}=\|v\|_{L^1(B)}+\|Dv\|_{L^2(B)}$ then the condition c.1) implies \beqno{vkbound}\|v_k\|_{L^1((0,1),W^{1,2}(B))}\le\int_0^1(\|v_k\|_{L^1(B)}+\|Dv_k\|_{L^2(B)})\,dt\le C(M).\eeq  

For any $h>0$, $t\in(-1,-h)$ and $k$ we denote $w_{k,h}(x,t)=v_k(x,t+h)-v_k(x,t)$. We thus have
\beqno{ww12} \int_{-1}^{-h}\|w_{k,h}\|_{W^{1,2}(B)}\,dt \le C(M).\eeq

Consider $p\in [1,2_*)$ and choose $l$ such that $lp'>N$. For any $\fg \in W^{l,p'}(B)$  $$\iidx{B}{\myprod{w_{k,h}, \fg}}=\iidx{B}{\myprod{\dspl{\int_{t}^{t+h}}(v_k)_t\,d\tau, \fg}}\le \dspl{\int_{t}^{t+h}}\iidx{B}{|(v_k)_t|}\, d\tau\|\fg \|_{L^\infty (B)}.$$ As $lp'>N$, by embedding theorems we have $\|\fg \|_{L^\infty (B)}\le C\|\fg \|_{W^{l,p'} (B)}$, the above implies $|\myprod{w_{k,h}, \fg}_{L^2(B)}|\le C\mu \|\fg \|_{W^{l,p'} (B)}$. This is to say \beqno{rsl}\|w_{k,h}\|_{W^{-l,p'}(B)}\le C\dspl{\int_{t}^{t+h}}\iidx{B}{|(v_k)_t|}\, d\tau.\eeq

Because $W^{1,2}(B)$ is compactly embedded in $L^p(B)$ and $L^p(B)$ is continuously embedded in $W^{-l,p'}(B)$, for any given $\mu>0$ we apply interpolating inequality to get
$$\|w_{k,h}\|_{L^p(B)}\le \mu \|w_{k,h}\|_{W^{1,2}(B)}+C(\mu)\|w_{k,h}\|_{L^{-l,p'}(B)}.$$
Raising the above to the power $q\ge1$, integrating over $t\in(-1,-h)$ and using \mref{ww12} and \mref{rsl}, we get
$$\barr{lll}\dspl{\int_{-1}^{-h}}{\|w_{k,h}\|_{L^p(B)}^q}\,dt&\le& \mu\dspl{\int_{-1}^{-h}}{\|w_{k,h}\|_{W^{1,2}(B)}^q}\,dt+C(\mu)\dspl{\int_{-1}^{-h}}{\|w_{k,h}\|_{W^{-l,p'}(B)}^q}\,dt\\&\le&\mu C(M)+C(\mu)\dspl{\int_{-1}^{-h}}\left[\dspl{\int_{t}^{t+h}}\iidx{B}{|(v_k)_t|}\,ds\right]^q\,dt.\earr$$
 From this, for any given $\eg >0$ we apply the continuity condition c.2) to the last integrand in the above (with $s=t$, $r=t+h$) to find $C(\eg)>0$ such that if $h<C(\eg)$ then
\beqno{wkheg}\dspl{\int_{-1}^{-h}}{\|w_{k,h}\|_{L^p(B)}^q}\,dt\le \eg.\eeq

We now see that for any $t_1,t_2\in(-1,0)$ the sequence $V_k(\cdot)=\int_{t_1}^{t_2}v(\cdot,s)ds$ is bounded in $W^{1,2}(B)$ so that it belongs to a fixed compact set in $L^p(B)$, thanks to \mref{vkbound}. Moreover, \mref{wkheg} clearly yields for all $k$ and $t\in(-1,0)$
$$\int_{-1}^{-h}\|v_k(\cdot,t+h)-v_k(\cdot,t)\|_{L^p(B)}^q\, dt \le O(h).$$
We then apply the well known compactness result of Simon (see \cite[Theorem 1]{Simon}) to see that $\{v_k\}$ is precompact in $L^q((-1,0),L^p(B))$. The lemma is proved. \eproof

 \brem{compactremm1}  The above lemma and its condition c.2) work well with  strong solutions  whose temporal derivatives are defined. Concerning weak solutions, we can replace c.2) by the following conditions which do not involve the derivatives $(v_k)_t$ (and then obtain a much better version of \cite[Lemma 3.3]{letrans}). We assume that there are sequences of functions $\{G_k\}$ and $\{f_k\}$ on $Q$ such that
 
\bdes
\item[c.2)]For all $\fg \in C^1_0(Q)$ there is a constant $C$  such that
$$\left|\itQ{Q}{\myprod{v_k,\fg_t }}\right|\le C\itQ{Q}{(|G_k||D\fg|+|f_k||\fg|)}.$$
\item[c.2')] For any given $\mu>0$ there is $C(\mu)$ such that if $-1<s<r<0$ and $r-s<C(\mu)$ then $$\int_s^r\iidx{B}{|G_k|}\, dt,\; \int_s^r\iidx{B}{|f_k|}\,dt\le \mu \quad \forall k.$$
 \edes

Let $l>N/p'+1$. By c.2) and because $C_0^1(B)\subset W^{l,p'}(B)$, the same argument in \cite[Lemma 3.2]{letrans}, with $u,G$ being $v_k,G_k$ and $f_k$ included, gives for any $-1<s<r<0$
$$\|v_k(\cdot,r)-v_k(\cdot,s)\|_{W^{-l,p'}(B)}\le C\int_s^r\iidx{B}{(|G_k|+|f_k|)}\,dt.$$
This is similar to \mref{rsl} in the proof and, together with the continuity condition c.2'), the proof can continue with this. \erem

\brem{compactrem} Concerning the continuity condition c.2) (or c.2' of \refrem{compactremm1}), we recall a well known result \cite[Corollary IV.11]{DunfordS} which shows that if a sequence $\{g_k\}$ converges {\em weakly} in $L^1(Q)$ then the following functions are absolutely continuous, {\em uniformly} in $k$.
$$A\to \itQ{A}{|g_k|},\; A\subset Q.$$ From this and the continuity of integrals, we then see that c.2) (respectively c.2')) is verified if the sequence $\{(v_k)_t\}$ (respectively $\{G_k\}, \{f_k\}$) converges weakly in $L^1(Q)$. 
In particular, if $\{(v_k)_t\}$ (respectively $\{G_k\}, \{f_k\}$) is a bounded sequence in $L^q(Q)$ for some $q>1$ then $\{(v_k)_t\}$ converges weakly in $L^1(Q)$ and c.2) is verified. In fact, this is the well known Aubin-Lions-Simon lemma (see \cite{Simon}).

 \erem
 
 In the proof we will frequently make use of the following interpolation Sobolev inequality 
 
 \blemm{Sobointineq} For any  $\eg>0$, $\bg\in(0,1]$ and $W\in W^{1,2}(\Og)$ we can find a constant $C(\eg,\bg)$  such that 
 \beqno{intineq}\|W\|_{L^q(\Og)}\le \eg\|DW\|_{L^p(\Og)} + C(\eg,\bg)\|W^\bg\|_{L^1(\Og)}^\frac{1}{\bg} \mbox{ for any $q\in[1,p_*)$}.\eeq \elemm 
 
\bproof By contradiction, assume that \mref{intineq} is not true then we can find $\eg_0>0$ and a sequence $\{W_n\}$ such that \beqno{intineqX}\|W_n\|_{L^q(\Og)}> \eg_0\|DW_n\|_{L^p(\Og)} + n\|W_n^\bg\|_{L^1(\Og)}^\frac{1}{\bg} \mbox{ for any $n$}.\eeq By scaling we can suppose that $\|W_n\|_{L^q(\Og)}=1$. The above implies that $\|DW_n\|_{L^p(\Og)}<1/\eg_0$ for all $n$. We see that $\{W_n\}$ is bounded in $W^{1,p}(\Og)$ so that, by compactness as $q<p_*$, we can assume that it converges to some $W$ in $L^q(\Og)$. Of course, $\|W\|_{L^q(\Og)}=1$. Meanwhile, \mref{intineqX} implies $\|W_n^\bg\|_{L^1(\Og)}\to0$ so that $\|W^\bg\|_{L^1(\Og)}=0$, this can be easily seen by H\"older's inequality and the H\"older continuity of the function $|x|^\bg $. Thus $W=0$ a.e on $\Og $ contradicting the fact that $\|W\|_{L^q(\Og)}=1$. The proof is complete. \eproof

\section{Existence of strong solutions}\eqnoset\label{res}

We discuss in this section the solvability of the following boundary and initial condition parabolic system.
\beqno{mainparaB-i}\left\{\barr{l} u_t-\Delta(P(u))= f(u),\mbox{ $(x,t)\in \Og\times(0,T_0)$},\\\mbox{$u=0$  on $\partial \Og\times(0,T_0)$},\\ u(x,0)=u_0(x) ,\quad x\in \Og. \earr\right.\eeq 

Firstly, we will apply the theory in \cite{dleANS,dlebook} to discuss the existence of strong solutions to this system when it is regular, i.e. PR) holds,
with initial data $u_{0}$ are in $W^{1,p_0}(\Og)$ for some $p_0>N$. We embed this system in the following family parameterized by $\sg\in[0,1]$
\beqno{mainparaBfam-i}\left\{\barr{l} u_t-\Delta(P(u))= \sg^2 f(u),\mbox{ $(x,t)\in \Og\times(0,T_0)$},\\\mbox{$u=0$ or $\frac{\partial u}{\partial \nu}=0$ on $\partial \Og\times(0,T_0)$},\\ u(x,0)=\sg u_{0}(x) ,\quad x\in \Og. \earr\right.\eeq

The existence of a strong solution to the regular system \mref{mainparaB-i} will be established under the crucial assumption that the {\em strong} solutions to the family \mref{mainparaBfam-i} {\em apriori} have small BMO norms (see \cite{FS, Gius}) in small balls (uniformly in $\sg\in[0,1]$). Namely, we consider the following property

\bdes\item[(Sbmo)] (Small BMO norm in small balls property) We say that a function $u:\Og\times(0,T_0)\to \RR^m$ satisfies (Sbmo) if for any given $\mu_0>0$ there is $R>0$ depending on the parameters in PR) and $\mu_0$ such that for any ball $B_R$ in $\RR^N$ with $\Og_R=B_R\cap \Og\ne\emptyset$ $$\sup_{t\in(0,T_0)}\|u(\cdot,t)\|_{BMO(\Og_R)} \le \mu_0.$$
\edes

Our first main result on the existence of strong solutions to the parabolic system \mref{mainparaB-i} is the following

\btheo{strongunithmB} Assume that PR), SG), and f) hold. Suppose further that {\em any strong} solution $u$ to the family \mref{mainparaBfam-i} apriori satisfy the following conditions. 
\bdes\item[a.1)] $\sg^{-1}u$ satisfies (Sbmo) {\em uniformly} in $\sg\in(0,1]$. \item[a.2)] There is a constant $C_0$ such that 
\beqno{llguucondaB-i} \sup_{t\in(0,T_0)}\|u(t)\|_{L^{1}(\Og)}\le C_0.\eeq 
\edes

Then there exists a strong solution $u$ to the system \mref{mainparaB-i}. 
\etheo

\bproof We apply \cite[Theorem 3.4.1]{dlebook} here by verifying its assumptions. First of all, we need to show that the number $\mathbf{\LLg}=\sup_{u\in\RR^m}\LLg(u)$,  with $\LLg(u)=|\llg_u(u)|/\llg(u)$, is finite.  Since $\llg(u)\ge \llg_{0}>0$, if $|u|$ is bounded then so is $\LLg(u)$. For large $|u|$ we use the assumption in P) that $|\llg_u(u)|\lesssim \llg(u)/|u|$ to see that $\LLg(u)\lesssim 1/|u|$ is also bounded. Hence, the number $\mathbf{\LLg}$ is finite.

Next, following \cite[Theorem 3.4.1]{dlebook}, we  consider the following family with $A(u)=P_u(u)$
\beqno{gensysfamB}\left\{\barr{l}u_t-\Div(A(\sg u)Du)=\sg f(\sg u),\quad (x,t)\in \Og\times(0,T_0),\\ u=0, \quad (x,t)\in \partial\Og\times(0,T_0),\\ u(x,0)=u_0(x)\quad x\in\Og.\earr\right.\eeq Multiplying $\sg>0$ to the equation in \mref{gensysfamB}, we see easily that $w=\sg u$ is a strong solution to 
\beqno{mainparaBfam}\left\{\barr{l} u_t-\Delta(P(u))= \sg^2 f(u),\mbox{ $(x,t)\in \Og\times(0,T_0)$},\\\mbox{$u=0$ on $\partial \Og\times(0,T_0)$},\\ u(x,0)=\sg u_{0}(x) ,\quad x\in \Og. \earr\right.\eeq First of all, the condition that strong solutions to \mref{gensysfamB} have small BMO norm in small balls of \cite[Theorem 3.4.1]{dlebook} is already assumed in a.1) that $u=\sg^{-1}w$ satisfies (Sbmo). We need only check the integrability conditions of \cite[Theorem 3.4.1]{dlebook}. From (Sbmo), we can take $\mu_0=1$ and find a fixed $R_1>0$ such that any strong solution $u$ to \mref{gensysfamB} satisfies
$$\sup_{t\in(0,T_0)}\|u(\cdot,t)\|_{BMO(B_{R_1})}\le 1.$$
For each $q\ge1$ and $t\in(0,T_0)$ it is well known (see \cite{Gius}) that $u(\cdot,t)$ is in $L^q(B_{R_1})$ and
$$\|u(\cdot,t)\|_{L^q(B_{R_1})}\le C(q,\|u(\cdot,t)\|_{BMO(B_{R_1})},\|u(\cdot,t)\|_{L^1(B_{R_1})}).$$

Since $\Og $ is bounded, by using a finite covering of finitely many balls of radius $R_1$ we deduce from the above   and the assumption \mref{llguucondaB-i} in a.2) of the theorem that for any $q\ge1$ there is a constant $C(q,R_1,C_0)$ which also depends  the geometry of $\Og $  such that
\beqno{BMOLq1}\sup_{t\in(0,C_0)}\|u(\cdot,t)\|_{L^q(\Og)}\le C(q,R_1,T_0).\eeq

From the polynomial growths of $\llg$ and $f$, we now see that $\llg (u)$, $|f(u)|\llg^{-1}(u)$ are in bounded by powers of $|u|$ so that their integrability conditions in \cite[Theorem 3.4.1]{dlebook} are verified by \mref{BMOLq1}. The last condition needs to be checked is 
\beqno{llgDu2bound1} \int_{0}^{T_0}\iidx{\Og}{|Du|^2}dt\le C_0(T_0)\eeq
for some constant $C_0(T_0)$. To prove this, we test the system \mref{gensysfamB} with $u$ and easily obtain
\beqno{llgDu2est}\int_0^{T_0}\iidx{\Og}{\llg(\sg u)|Du|^2}dt\le \iidx{\Og}{\sg |u_0|^2}+ \int_0^{T_0}\iidx{\Og}{\sg \myprod{f(\sg u),u}}dt.\eeq  By the polynomial growth of $f(u)$, the integrand on right hand side of the above is bounded by a polynomial in $|u|$. By \mref{BMOLq1}, we conclude that the right hand side of \mref{llgDu2est} is bounded uniformly in $\sg\in[0,1]$.
On the other hand, as $\llg(\sg u)$ is bounded from below by $\llg_0>0$, we obtain \mref{llgDu2bound1}. The  proof is complete. \eproof

\section{Existence of weak solutions} \eqnoset\label{weaksol}
Next,  we study the existence of a weak solution to the  following boundary and initial condition problem.
\beqno{mainpara}\left\{\barr{l} u_t-\Delta(P(u))=f(u),\mbox{ $(x,t)\in \Og\times(0,T_0)$},\\\mbox{$u=0$ on $\partial \Og\times(0,T_0)$},\\ u(x,0)=u_{0}(x) ,\quad x\in \Og, \earr\right.\eeq 
where $P$ is only assumed to satisfy the condition P), i.e. $\llg_0$ can be 0. 

We state the standard definition of weak solutions here.

{\bf Definition W):} We say that $u$ is a weak solution to \mref{mainpara} in $Q=\Og\times(0,T_0)$ if $u\in L_{loc}^1(Q)$ and $P(u)\in L_{loc}^1(0,T_0:W^{1,1}(\Og))$; and
for any $\eta\in C^1(\bar{Q})$, $\eta=0$ on $\partial\Og\times (0,T)$ and $\Og \times\{T\}$ the following holds 
$$\itQ{Q}{(-\myprod{u,\eta_t}+\myprod{DP(u),D\eta})}=\iidx{\Og}{u_{0}\eta(x,0)}+\itQ{Q}{\myprod{f(u),\eta}}.$$

Inspired by the scalar porous media model, we assume further that \bdes \item[Ph)] $P^{-1}$ exists and is H\"older continuous for some $\ag_P\in(0,1]$: There is a constant $[P]_{\ag_P} >0$ such that $|P^{-1}(u)-P^{-1}(v)|\le [P]_{\ag_P}|u-v|^{\ag_P}$ for all $u,v\in \RR^m$. Equivalently, \beqno{Phequiv}|u-v|\le [P]_{\ag_P}|P(u)-P(v)|^{\ag_P}\mbox{ for all $u,v\in \RR^m$}.\eeq

\edes

An example of such $P$ is $P(u)=|u|^ku$ for some $k>0$. Then $P^{-1}(u)=|u|^\frac{-k}{1+k}u$ which is H\"older continuous with the exponent $\ag_P=1/(k+1)$, this is the porous media model we discussed in the Introduction. The map $P$ defined for the generalized SKT) system in the Introduction also satisfies this condition. Indeed, away from the singular point $u=0$, $P$ is Lipschitz because $P_u^{-1}$ exists and bounded. At $u=0$, it is clear that \mref{Phequiv} holds because $|P(v)|\ge C|v|^{k+1}$ for some positive constant $C$.

We will obtain a weak solution to the degenerate system as the limit of a sequence of strong solutions to regularized systems. To this end, let $\{\llg_{0,n}\}$ be a sequence in $(0,1)$ and $\lim_{n\to\infty}\llg_{0,n}=0$. We denote $P_{n}(u)=\llg_{0,n} u + P(u)$   and consider the following approximation systems with initial data $u_{0,n}$ being in $W^{1,p_0}(\Og)$ for some $p_0>N$.
\beqno{regpara}\left\{\barr{l} u_t-\Delta(P_{n}(u))=f(u),\quad (x,t)\in \Og\times(0,T_0),\\\mbox{$u=0$ on $\partial \Og\times(0,T_0)$},\\ u(x,0)=u_{0,n}(x) ,\quad x\in \Og. \earr\right.\eeq 

The system \mref{regpara} satisfies PR) because $\llg_{0,n}>0$. Following \reftheo{strongunithmB},
for each $n$ we embed \mref{regpara} in the following family of systems parameterized by $\sg\in[0,1]$
\beqno{parafamequiv}\left\{\barr{l} u_t-\Delta(\llg_{0,n} u + P(u))=\sg^2 f(u),\quad (x,t)\in \Og\times(0,T_0),\\\mbox{$u=0$ on $\partial \Og\times(0,T_0)$},\\ u(x,0)=\sg u_{0,n}(x) ,\quad x\in \Og. \earr\right.\eeq

If strong solutions to the above system apriori satisfy the assumption a.1) and a.2) of \reftheo{strongunithmB} then
we obtain a sequence of strong solutions $\{u_n\}$ for \mref{regpara}. However, in order to pass to the limit to obtain the existence of a weak solution to our degenerate system \mref{mainpara}, we have to assume that these strong solutions satisfy a bit stronger integrability condition than a.2) of \reftheo{strongunithmB} {\em uniformly in $n$} (see \mref{llguuconda-i} below).

Concerning the initial condition of \mref{mainpara} and \mref{regpara}, we also assume that \bdes\item[IC)] There exists a sequence $\{u_{0,n}\}$ in $C^1(\Og)$ which converges to $u_0$ in $L^1(\Og)$. Furthermore, there is a constant $C_0$ such that for all $n$ \beqno{llguDu0conda-i}\|(\llg_{0,n}+\llg(u_{0,n}))Du_{0,n}\|_{L^{2}(\Og)}\le C_0.\eeq\edes

\btheo{strongweakthm} Assume P), Ph), IC), f). Let $\{\llg_{0,n}\}$ be a sequence in $(0,1)$ and $\lim_{n\to\infty}\llg_{0,n}=0$. Consider the family \mref{parafamequiv} and assume that its strong solutions apriori satisfy the condition a.1) of \reftheo{strongunithmB} {\bf for each $n$}. 

Assume also that  there is a constant $q_0>N/2$ and $C_1$ such that any strong solutions $u_n$ of \mref{regpara} satisfy

\beqno{llguuconda-i}\sup_{t\in(0,T_0)}\|\llg(u_n)\|_{L^{q_0}(\Og)},\; \sup_{t\in(0,T_0)}\|u_n\|_{L^{1}(\Og)}\le C_1.\eeq 

Then there exists a weak solution $u$ to the system \mref{mainpara}. \etheo

We should emphasize that the condition a.1) of \reftheo{strongunithmB} on the property (Sbmo) is assumed for each $n$ in order to obtain the strong solutions to the regular systems \mref{regpara} and this condition  is uniform only in $\sg\in(0,1]$ but not in $n$. Meanwhile, the integrability condition \mref{llguuconda-i} is assumed to be uniform in $n$.

In order to pass to the limit to obtain the existence of a weak solution to our degenerate system \mref{mainpara}, we have to to establish  {\em uniform} estimates for these strong solutions $u_n$ under the integrability conditions \mref{llguDu0conda-i} and \mref{llguuconda-i} of \reftheo{strongweakthm}.

The following proposition provides the needed uniform estimates.

\bprop{strongunithm} Assume P) and f). 
Assume also that  there are constants $q_0>N/2$ and $C_0,C_1$ such that the initial data $u_{0,n}\in C^1(\Og)$ and the corresponding strong solutions $u_n$ of \mref{regpara} satisfy

\beqno{llguDu0conda}\|u_{0,n}\|_{L^{\infty}(\Og)},\;\|\llg(u_{0,n})Du_{0,n}\|_{L^{2}(\Og)}\le C_0,\eeq
\beqno{llguuconda}\sup_{t\in(0,T_0)}\|\llg(u_n)\|_{L^{q_0}(\Og)},\; \sup_{t\in(0,T_0)}\|u_n\|_{L^{1}(\Og)}\le C_1.\eeq 

Then  there are constants $C(C_0,C_1)$ and $q_1>1$ such that for every $n$
\beqno{u2est0}\|u_n(t)\|_{L^2(\Og)} \le C(T_0,C_1),\eeq
\beqno{Du-uniboundpara} \sup_{t\in[0,T_0]}\iidx{\Og}{(\llg_{0,n}^2+\llg^2(u_n))|Du_n|^2} \le C(C_0,C_1),\eeq
\beqno{ut-uniboundpara}  \itQ{\Og\times[0,T_0]}{\llg(u_n)|(u_n)_t|^2}\le C(C_0,C_1),\eeq and  \beqno{festa}\sup_{t\in(0,T_0)}\iidx{\Og}{|f(u_n)|^{q_1}}\le C(C_0,C_1).\eeq
\eprop

The uniform estimates \mref{Du-uniboundpara}-\mref{festa} will come from following lemmas which discuss the estimates for strong solutions of \beqno{genpara}\left\{\barr{l} u_t-\Delta(\ccP(u))=f(u),\quad (x,t)\in \Og\times(0,T_0),\\u=0\mbox{ on $\partial\Og \times(0,T_0)$}. \earr\right.\eeq Here, $\ccP=[\ccP_i]_{i=1}^m$ is a $C^2$ map on $\RR^m$ and satisfies the condition P).
In the lemmas and their proof, we will denote $\ccA(u)=\ccP_u(u)$ and the ellipticity function $\llg$ for $\ccA$ by $\llg_\ccA$.

To begin, note that the ellipticity condition \mref{A1} of P) and Young's inequality imply $$\llg_\ccA(u)|Du|^2\le \myprod{\ccA(u)Du,Du}=\myprod{D\ccP(u),Du}\le \frac12\llg_\ccA^{-1}(u)|D\ccP(u)|^2+\frac12\llg_\ccA (u)|Du|^2.$$ We then have $\llg_\ccA(u)|Du|^2\le \llg_\ccA^{-1}(u)|D\ccP(u)|^2$ so that $\llg_\ccA(u)|Du|\le |D\ccP(u)|$. Of course, $|D\ccP(u)|= |\ccA(u)Du|\le C_*\llg_\ccA(u)|Du|$. Hence,
\beqno{DuleDU} \llg_\ccA(u)|Du|\sim |D\ccP(u)|.\eeq

The first lemma provides a differential (or Gronwall) inequality for $\|\ccA(u)Du\|_{L^2(\Og)}$.
\blemm{ldulemm0z} Let $u$ be a strong solution to \mref{genpara}. For any $t\in(0,T_0)$ \beqno{Aeqn}\iidx{\Og\times\{t\}}{\llg_\ccA(u)|u_t|^2} +\frac{d }{d t}\iidx{\Og\times\{t\}}{|\ccA(u)Du|^2}\le C\iidx{\Og\times\{t\}}{\llg_\ccA(u)|f(u)|^2}.\eeq

\elemm

\bproof
Because $u$ is a strong solution, we can test the system with $\ccP(u)_t$. This means we multiply the $i^{th}$ equation of the system by $(\ccP_i(u))_t$  and integrate by parts in $x$ over $\Og$. Summing the results, we get for any $t\in(0,T_0)$
$$\iidx{\Og\times\{t\}}{(\myprod{\ccP(u)_t,u_t} +\myprod{D(\ccP(u)),D(\ccP(u)_t)})}=\iidx{\Og\times\{t\}}{\myprod{f(u),\ccP(u)_t}}.$$

As $D(\ccP(u)_t)=(D\ccP(u))_t$,  we have $\myprod{D(\ccP(u)),D(\ccP(u)_t)}=\frac12\frac{\partial }{\partial t}(|D(\ccP(u))|^2)$ so that
\beqno{ztemp}\iidx{\Og\times\{t\}}{[\myprod{\ccA(u)u_t,u_t} +\frac12\frac{\partial }{\partial t}(|D(\ccP(u))|^2)]}=\iidx{\Og\times\{t\}}{\myprod{f(u),\ccA(u)u_t}}.\eeq 
We now use the ellipticity of $\ccA(u)$ in the first integrand on the left hand side of \mref{ztemp} to have $\myprod{\ccA(u)u_t,u_t}\ge\llg_\ccA(u)|u_t|^2$. Also, as  $|\ccA(u)|\le C\llg_\ccA(u)$, we use Young's inequality to find a constant $C(\eg)$ such that for any $\eg>0$ we can estimate the second integrand on the right hand side as follows $|\myprod{f(u),\ccA(u)u_t}|\le  \eg\llg_\ccA(u)|u_t|^2 + C(\eg)\llg_\ccA(u)|f(u)|^2$.
Using these facts in \mref{ztemp} with  sufficiently small $\eg$, we get \mref{Aeqn}.  \eproof

In order to estimate the integral of $\llg_\ccA(u)|f(u)|^2$ in \mref{Aeqn} we need the following lemma.

\blemm{llgfulemma} Assume that $|u||(\llg_\ccA)_u(u)|\lesssim \llg_\ccA(u)$ and  that there are constants  $q_0>N/2$ and $C_1$ such that 
\beqno{llguucond}\|\llg_\ccA(u)\|_{L^{q_0}(\Og)},\; \|u\|_{L^{1}(\Og)}\le C_1.\eeq  

There is $q\in(2,2_*)$ such that for any given $\eg $ there is a constant $C(\eg ,C_1)$ such that \beqno{llguuz}\left(\iidx{\Og}{(\llg_\ccA(u)|u|)^{q}}\right)^\frac{2}{q}\le \eg \iidx{\Og}{|\llg_\ccA(u)Du|^2}+C(\eg, C_1),\eeq
and
\beqno{llgfest1}\iidx{\Og}{|u|^2\llg_\ccA^3(u)}\le \eg \iidx{\Og}{|\llg_\ccA(u)Du|^2}+C(\eg ,C_1).\eeq
\elemm

\bproof  From the asumptions on $q_0$ 
it is clear that we can find $q\in(2,2_*)$ such that $\frac{N}{2}<(\frac{q}{2})'=q/(q-2)\le q_0$. By the H\"older inequalitie and the assumption \mref{llguucond}, we have  \beqno{llg3u2}\barr{lll}\iidx{\Og}{\llg_\ccA^3(u)|u|^2}&\le&\left(\iidx{\Og}{(\llg_\ccA(u)|u|)^{q}}\right)^\frac{2}{q}\|\llg_\ccA(u)\|_{L^{(\frac{q}{2})'}(\Og)}\\&\le& C_1\left(\iidx{\Og}{(\llg_\ccA(u)|u|)^{q}}\right)^\frac{2}{q},\earr\eeq

Thus, \mref{llgfest1} folllows from  \mref{llguuz}, which we will prove below.

Because $q<2_*$, we can apply the interpolation inequality \mref{intineq} to estimate the integral of $(\llg_\ccA(u)|u|)^q$. First of all, we note that $|D(\llg_\ccA(u)|u|)|\lesssim  \llg_\ccA(u)|Du|$ thanks to  the assumption $|u||(\llg_\ccA)_u(u)|\lesssim \llg_\ccA(u)$ of the lemma. We then have, by the interpolation inequality, for any given $\eg,\bg>0$
\beqno{llguu}\left(\iidx{\Og}{(\llg_\ccA(u)|u|)^{q}}\right)^\frac{2}{q}\le \eg\iidx{\Og}{|\llg_\ccA(u)Du|^2}+C(\eg,\bg)\left(\iidx{\Og}{(\llg_\ccA(u)|u|)^{\bg}}\right)^\frac{2}{\bg}.\eeq

As $(\llg_\ccA(u)|u|)^{\bg}= (\llg_\ccA^3(u)|u|^2)^{\bg/3}|u|^{\bg/3}$, we have (if $\bg<3$) $$\iidx{\Og}{(\llg_\ccA(u)|u|)^{\bg}}\le \left(\iidx{\Og}{\llg_\ccA^3(u)|u|^{2}}\right)^\frac{\bg}{3}\left(\iidx{\Og}{|u|^\frac{\bg}{3-\bg}}\right)^\frac{3-\bg}{3}.$$ We choose $\bg<3/2$   so that $\bg/(3-\bg)\le 1$. By the assumption  \mref{llguucond} and Young's inequality, we obtain from the above that for any $\eg_0>0$ $$\left(\iidx{\Og}{(\llg_\ccA(u)|u|)^{\bg}}\right)^\frac{2}{\bg}\le C_1\left(\iidx{\Og}{\llg_\ccA^3(u)|u|^{2}}\right)^\frac{2}{3}\le \eg_0 \iidx{\Og}{\llg_\ccA^3(u)|u|^{2}}+C(\eg_0,C_1).$$ We estimate the last integral by \mref{llg3u2} and then use  \mref{llguu} to get
$$\barr{lll}\left(\iidx{\Og}{(\llg_\ccA(u)|u|)^{\bg}}\right)^\frac{2}{\bg}&\le& C_1\eg_0\eg\iidx{\Og}{|\llg_\ccA(u)Du|^2}+\\&&C_1\eg_0C(\eg,\bg)\left(\iidx{\Og}{(\llg_\ccA(u)|u|)^{\bg}}\right)^\frac{2}{\bg}
+C(\eg,C_1).\earr$$ Clearly, for any given $\eg >0$ we can find $\eg_0$ such that   $C_1\eg_0\max\{1,C(\eg ,\bg)\}<1/2$, the second integral on the right hand side can then be absorbed into the left. We then have $$\left(\iidx{\Og}{(\llg_\ccA(u)|u|)^{\bg}}\right)^\frac{2}{\bg}\le \eg \iidx{\Og}{|\llg_\ccA(u)Du|^2}+C(\eg ,C_1).$$

Using this in \mref{llguu}, we obtain \mref{llguuz} and complete the proof. 
\eproof

\blemm{llgfulemma1} Assume as in \reflemm{llgfulemma}. We also find a constant $C(C_1)$ such that
\beqno{u2est}\|u(t)\|_{L^2(\Og)} \le C(T_0,C_1),\eeq
\beqno{llgfest}\iidx{\Og}{|u|^2\llg_\ccA(u)}\le C\iidx{\Og}{|\llg_\ccA(u)Du|^2}+C(C_1).\eeq
\elemm

\bproof We test the system for $u$ and easily obtain for any $t\in(0,T_0)$ that
$$\iidx{\Og\times\{t\}}{|u|^2}+\itQ{\Og\times(0,t)}{\llg_\ccA|Du|^2}\le \iidx{\Og}{|u_0|^2}+\itQ{\Og\times(0,t)}{|f(u)||u|}.$$

As $|f(u)||u|\le\frac12(|f(u)|^2+|u|^2)$ and $|f(u)|^2\le |u|^2+|\llg_\ccA|^2|u|^2$, we can make use of \mref{llguuz}, with sufficiently small $\eg $, to arrive at
$$\iidx{\Og\times\{t\}}{|u|^2}\le \iidx{\Og}{|u_0|^2}+C(C_1)+C\itQ{\Og\times(0,t)}{|u|^2}.$$ This is a Gronwall inequality for $\|u(t)\|_{L^2(\Og)}^2$ and it yields \mref{u2est}. 

On the other hand, by H\"older's inequality
\beqno{llgu2}\iidx{\Og}{\llg_\ccA(u)|u|^2}\le\left(\iidx{\Og}{(\llg_\ccA(u)|u|)^{2}}\right)^\frac{1}{2}\left(\iidx{\Og}{|u|^{2}}\right)^\frac{1}{2}.\eeq Combining this with \mref{llguuz} and \mref{u2est}, we obtain \mref{llgfest}. The proof is complete.
\eproof

We are now ready to prove  \refprop{strongunithm}.

{\bf Proof of \refprop{strongunithm}:}   Let $\ccA(u)=\llg_{0,n}I+A(u)$. The estimate \mref{u2est0} for $\|u_n(t)\|_{L^2(\Og)}$ comes from \mref{u2est}. The integrability condition \mref{llguuconda} implies \mref{llguucond} of \reflemm{llgfulemma}  so that the estimates \mref{llgfest1}, \mref{llgfest} for the integrals of $|u|^2\llg_\ccA^3(u)$ and $|u|^2\llg_\ccA(u)$ hold. The growth condition \mref{fumore} gives $\llg_\ccA(u)|f(u)|^2\lesssim  |u|^2\llg_\ccA(u)+|u|^2\llg_\ccA^3(u)$ so that $$\iidx{\Og}{\llg_\ccA(u)|f(u)|^2}\le C\iidx{\Og}{|\ccA(u)Du|^2}+C(C_1).$$
Here, we used the fact that $|\llg_\ccA(u)Du|\sim |\ccA(u)Du|$ (see \mref{DuleDU}).
We can use the above in \mref{Aeqn} of \reflemm{ldulemm0z} to obtain
\beqno{Aeqnz}\iidx{\Og\times\{t\}}{\llg_\ccA(u)|u_t|^2} +\frac{d }{d t}\iidx{\Og\times\{t\}}{|\ccA(u)Du|^2}\le C\iidx{\Og\times\{t\}}{|\ccA(u)Du|^2} +C(C_1).\eeq

Define $y(t):=\|\ccA(u)Du\|_{L^2(\Og\times\{t\})}^2$. We obtain from \mref{Aeqnz} that $y'\le Cy+C(C_1)$. By Gronwall's inequatity, we see that any strong solution $u_n$ of \mref{regpara} satisfies $$\|\ccA(u)Du\|_{L^2(\Og\times\{t\})}^2\le C(\|Du_{0,n}\|_{L^2(\Og)},\|u_{0,n}\|_{L^\infty(\Og)})+C(C_1) \mbox{ for any $t\in(0,T_0)$}.$$ 

As  $\llg_\ccA(u)=\llg_{0,n}+\llg(u)$ and $|\ccA(u)Du|\lesssim \llg_\ccA(u)|Du|$, we use the above and the assumptions on the initial condition \mref{llguDu0conda} on the initial data $u_{0,n}$ to prove \mref{Du-uniboundpara} of the proposition. 

Next, by integrating \mref{Aeqnz}, we then have for all $t\in(0,T_0)$ that
$$ \barr{ll}\lefteqn{\int_t^{T_0}\iidx{\Og}{\llg_\ccA(u)|u_t|^2}\,ds+\iidx{\Og\times\{T_0\}}{|\ccA(u)Du|^2}\le}\hspace{2cm}&\\&  \iidx{\Og\times\{t\}}{|\ccA(u)Du|^2}+C(C_1)\le C(C_0,C_1).\earr
$$
Letting $t\to0$ and using the assumption \mref{llguDu0conda} (and \refrem{limDuat0} after this proof) on the initial data, we obtain $$ \int_0^{T_0}\iidx{\Og}{\llg_\ccA(u)|u_t|^2}\,ds\le C(C_0,C_1),$$ 
and prove \mref{ut-uniboundpara}. 

Finally, let $q_1=\min\{q,2\}>1$ with $q$ being the exponent in \mref{llguuz}. By \mref{fumore}
\beqno{fest}\iidx{\Og}{|f(u)|^{q_1}}\le C\iidx{\Og}{(|u|^{q_1}+|u|^{q_1}\llg^{q_1}(u))},\eeq
so that the estimate \mref{festa} for $f(u)$ comes from the bound \mref{u2est0} and the inequality \mref{llguuz} of \reflemm{llgfulemma} in combination with the bound \mref{Du-uniboundpara}. The proof is complete. \eproof

We are now ready to present the proof of \reftheo{strongweakthm} on the existence of a weak solution to the degenerate systems.

{\bf Proof of \reftheo{strongweakthm}:}  Consider the sequence of strong solutions $\{u_n\}$ obtained from \reftheo{strongunithmB}, with initial data $u_{0,n}$. This sequence exists because we are assuming the conditions a.1) and a.2) for each $n$ here and \reftheo{strongunithmB} applies.

Denote $U_n:=P(u_n)$. For any $q\in(1,2)$, because $$|(U_n)_t|^q\lesssim \llg_n(u_n)^q|(u_n)_t|^q=\llg_n(u_n)^\frac q2\llg_n(u_n)^\frac q2|(u_n)_t|^q,$$ we can apply H\"older's inequality  to obtain  for $Q=\Og\times[0,T_0]$ that
$$ \itQ{Q}{|(U_n)_t|^q}\le \left(\itQ{Q}{\llg_n(u_n)^\frac{q}{2-q}}\right)^{1-\frac q2} \left(\itQ{Q}{\llg_n(u_n)|(u_n)_t|^2}\right)^\frac q2.$$
As we assume that $\|\llg(u_n)\|_{L^{q_0}(\Og)}$ is uniformly bounded for some $q_0>N/2\ge1$, there is $q>1$ such that  $\frac{q}{2-q}\in (1,q_0)$ and therefore the first integral on the right hand side is bounded uniformly by a constant. By \mref{ut-uniboundpara} of \refprop{strongunithm}, the second integral is also bounded. Thus, $\{(U_n)_t\}$ is bounded in $L^q(Q)$ and
we can use \reflemm{compactlemm} to see that $\{U_n\}$ is precompact in $L^p([0,T_0],L^{p}(\Og))$ for any given $p\in (1,2_*)$.

Hence, for $p=2$ we can find a subsequence of $\{U_n\}$ such that, after relabeling
\beqno{Unconv} U_n \to U \quad \mbox{ in } X:=L^2([0,T_0],L^{2}(\Og)).\eeq
Via a subsequence again, we can assume that $U_n(t)\to U(t)$ in $L^2(\Og)$ a.e in $[0,T_0]$. By \mref{u2est0} $\{u_n\}$ is bounded in $X$ so that we can also assume that it converges weakly to some $u(t)\in X$.  Since $P^{-1}$ exists, using the uniqueness of weak limits, we have $u(t):=P^{-1}(U(t))$. In fact, using the H\"older continuity of $P^{-1}$ in Ph) and \mref{Unconv} we easily see that $u_n\to u$ in $L^q(Q)$ for any $q\le 2/\ag_P$, where $\ag_P$ is the H\"older exponent of $P^{-1}$.

Furthermore, \mref{Du-uniboundpara} shows that the sequence $\{D(\llg_{0,n} u_n(t)+U_n(t))\}$ is bounded in $L^2(\Og)$ so that it converges weakly. Note that for any $\fg \in C_0^\infty(\Og)$ as $n\to\infty$
$$\iidx{\Og}{\myprod{DP(u_n(t)),\fg}}=-\iidx{\Og}{\myprod{P(u_n(t)),D\fg}}\to -\iidx{\Og}{\myprod{P(u(t)),D\fg}}.$$
Thus, $D(\llg_{0,n} u_n+U_n)$ converges weakly to $DP(u)$ in the sense of distribution. In fact, the bound in \mref{Du-uniboundpara} for $\{D(\llg_{0,n} u_n(t)+U_n(t))\}$ in $L^2(\Og)$ and density show that $D(\llg_{0,n} u_n+U_n)$ converges weakly to $DP(u)$ in $L^2(\Og)$.

On the other hand, by \mref{festa}, $f(u_n)$ is unformly bounded in $L^{q_1}(\Og)$ for some $q_1>1$, it converges weakly to $f(u)$ in $L^{q_1'}(\Og)$ (see also \refrem{fconvrem} below).

For any $\eta\in C^1(\bar{Q})$, $\eta=0$ on $\partial\Og\times (0,T)$ and $\Og \times\{T\}$,  we multiply $\eta$ to the equation of the strong solution $u_n$ and derive
$$\barr{ll}\lefteqn{\itQ{Q}{(-\myprod{u_n,\eta_t}+\myprod{D(\llg_{0,n} u_n+U_n),D\eta})}=}\hspace{3cm}&\\&\iidx{\Og}{u_{0,n}\eta(x,0)}+\itQ{Q}{\myprod{f(u_n),\eta}}.\earr$$ 

Let $n\to\infty$. By the convergences established above and the condition on the initial data in IC) we obtain 
$$\itQ{Q}{(-\myprod{u,\eta_t}+\myprod{DP(u),D\eta})}=\iidx{\Og}{u_{0}\eta(x,0)}+\itQ{Q}{\myprod{f(u),\eta}}.$$

We see that $u$ is a weak solution. The proof is complete. \eproof

\brem{fconvrem} As we considering polynomial growth data in this paper, we can assume that $|P(u)|\sim |u|^{k+1}$ for some $k>0$.  So that the ellipticity function $\llg (u)\lesssim |u|^k+1$. It has been seen from the proof that  $u_n\to u$ in $L^{q}(Q)$ for some $q>2(k+1)$ because $P(u_n)$ converges in $L^2(Q)$ and $\ag_P<1$. By the Riesz-Fisher theorem we can extract a subsequence of $u_n$ and assume that there is a function $\hat{u}\in L^{q}(Q)$ such that $u_n\to u$ a.e. in $Q$ and $|u_n|\le \hat{u}$ for all $n$. Thus, by f), we have $|f(u_n)|\le |u_n|+|u_n|\llg (u_n))\lesssim \hat{u}^{k+1}+\hat{u}$, a function in $L^2(Q)$. Because $f$ is continuous, we have $f(u_n)\to f(u)$ a.e. in $Q$.  By Dominated convergence theorem, we see that $f(u_n)\to f(u)$ in $L^2(Q)$. \erem

\section{The planar case $N=2$:}\eqnoset\label{planar}
The crucial condition (Sbmo) in a.1) of \reftheo{strongunithmB} must be established in order to establish the existence of a sequence of strong solutions to the approximation systems.  This condition is not easy to validate in general. However, when $N=2$,  \refprop{strongunithm} provides a bound for $\sup_{t\in(0,T_0)}\|Du_n\|_{L^{2}(\Og)}$ and allows us
to verify the (Sbmo) property under a very weak a priori integrability condition of strong solutions. On the other hand, as the  the H\"older continuity of the strong solutions $u_n$ obtained by \reftheo{strongunithmB} is not uniform when $\llg_{0,n}\to0$  so that this regularity cannot pass to that of the weak solution $u$ found in \reftheo{strongweakthm}. 
At least, we can show that this weak solution $u$  is VMO (Vanishing Mean Oscillation). That is, $$ \limsup_{R\to0} \|u\|_{BM0(\Og_R(x,t))}=0, \quad \forall (x,t)\in\Og_{R}\times(0,T_0).$$

\btheo{Nis2theo} Let $N=2$. Assume that P), Ph), f) and \mref{fumore} hold. For any given $\llg_{0,n}>0$ assume that there are constants $q_0>1$ and $C_1$ such that strong solutions of \mref{parafamequiv}
apriori satisfy
\beqno{llguucondaN2}\|\llg(u)\|_{L^{q_0}(\Og)},\; \|u\|_{L^{1}(\Og)}\le C_1.\eeq

Then there exists a weak solution $u$ to \mref{mainpara}. Moreover, $u$ is VMO.
\etheo

\bproof First of all, we show that the condition a.1) of \reftheo{strongunithmB} holds so that strong solutions of \mref{parafamequiv} exist for $\sg=1$. Consider a strong solution $u$ of the family  \mref{parafamequiv}, $\sg\in(0,1]$. Under the condition \mref{llguucondaN2} \refprop{strongunithm} applies here with $q_0>1$ (because $N=2$) and $f(u)$ being $\sg^2 f(u)$. We then obtain from \mref{Du-uniboundpara} $$\sup_{t\in[0,T_0]}\iidx{\Og}{(\llg_{0,n}^2+\llg^2(u))|Du|^2} \le \sg^2C(C_1)$$ and this implies $$\sup_{t\in[0,T_0]}\iidx{\Og}{|D(\sg^{-1}u)|^2} \le \llg_{0,n}^{-2}C(C_1).$$ As $N=2$, a simple use of Poincar\'e's inequality, the continuity of integral  and the last estimate show that $\sg^{-1}u$ satisfies the (Sbmo) condition (uniformly in $\sg\in(0,1]$). Thus, a.1) is verified. The condition a.2) of \reftheo{strongunithmB} is assumed in \mref{llguucondaN2} here. We obtain a sequence of strong solutions $\{u_n\}$ to \mref{regpara} for $\sg=1$. Using \reftheo{strongweakthm} and letting $\llg_{0,n}\to0$ we then obtain a weak solution $u$. 

To finish the proof we will need only show that $u$ is VMO. First of all,
by \mref{Du-uniboundpara} and because $|D(P(u_n))|\lesssim \llg(u_n)|Du_n|$, the strong solutions satisfy $$\sup_{(0,T_0)}\|D(P(u_n))\|_{L^2(\Og)} \le C(C_1).$$ 

Let $U_n=P(u_n)$. For any $q>1$, from the minimizing property of average, it is well known that there is a constant $c(q)$ such that
$$\iidx{\Og_R}{|u_n-(u_n)_R|^q} \le c(q)\iidx{\Og_R}{|u_n-P^{-1}(U_n)_R|^q}.$$ 
We use the H\"older property  of $P^{-1}$ in Ph) to estimate the last integral. 
\beqno{unUn}\iidx{\Og_R}{|P^{-1}(U_n)-P^{-1}(U_n)_R|^q}\le [P^{-1}]_{\ag_P}^q\iidx{\Og_R}{|U_n-(U_n)_R|^{q\ag_P}}.\eeq

By the Poincar\'e-Sobolev inequality and the uniform continuity of the integrals,  for any $\mu_0>0$, there is $R>0$ depends only on $\mu_0$ such that
$$ R^{-2}\iidx{\Og_R}{|U_n-(U_n)_R|^{2}}\le \iidx{\Og_R}{|DU_n|^{2}}\le \mu_0,\quad \forall n.$$

Take $q=2/\ag_P$. We combine the above estimates to obtain 
$$ R^{-2}\iidx{\Og_R}{|u_n-(u_n)_R|^q}\le C(\ag_P,[P^{-1}]_{\ag_P})\mu_0,\quad \forall n.$$

From the proof of \reftheo{strongweakthm},  $U_n(t)\to U(t)$ in $L^2(\Og)$ for a.e. $t\in(0,T_0)$ so that $u_n\to u$ in $L^q(\Og)$, because of \mref{unUn}.  Letting $n\to\infty$ in the above estimate, we see that $u$ satisfies it too. By the equivalence of BMO norm definitions, we have $[u]_{BMO(\Og_R)}\le C(\ag_P,[P^{-1}]_{\ag_P})\mu_0$. As $\mu_0$ can be arbitrarily small, if $R$ is, $u$ is VMO. The proof is complete. \eproof

We are now ready to provide the existence part of a weak solution to the degenerate \mref{iSKT} stated in the Introduction. We just need to establish the bounds of the norms in \mref{llguucondaN2} for strong solutions to the nondegenerate family \mref{parafamequiv}.

Under the condition \mref{agcondyagi} on $\ag_{ij}$'s, namely $ \ag_{21}^2<8\ag_{11}\ag_{12}$ and $\ag_{12}^2<8\ag_{22}\ag_{21}$, Yagi showed in \cite{yag} that if the initial data $u_0=[u_1(x,0), u_2(x,0)]^T$ are nonnegative then the strong nonegative solution to the nondegenerate  system \mref{e0} are also nonnegative and there is $\llg(u)\sim |u|$ such that for some positive constant $c_\ag$ depending on $\ag_{ij}$'s \beqno{iSKTP_u} \myprod{P_u Du,Du}\ge c_\ag|u||Du|^2, \quad \mbox{where } P(u):=\left[\barr{c}u_1(\ag_{11}u_1+\ag_{12}u_2)\\
u_2(\ag_{21}u_1+\ag_{22}u_2)\earr\right].\eeq

We just need to show that $\|u\|_{L^2(\Og)}$ is bounded uniformly with respect to $\sg$ to establish \mref{llguucondaN2} (for $q_0=2$). This is exacly what will be done in the following two lemmas.

We consider first the competitive (SKT).
\blemm{CufuSKTlem} Assume \beqno{SKTcompete} \myprod{f(u),u} \le C|u|^2,\;\forall u=[u_1,u_2]^T\in \RR^2, u_1,u_2\ge 0.\eeq 

Then $\|u(t)\|_{L^2(\Og)}^2$ is bounded by a constant $c(T_0,\|u_0\|_{L^2(\Og)})$. \elemm

\bproof 
We test the system with $u$ and use \mref{iSKTP_u}, \mref{SKTcompete} to easily get for any $T\in(0,T_0)$,  $\sg\in(0,1)$ and $Q=\Og\times(0,T)$, dropping the integral of $|Du|^2$ on the left
\beqno{keyestSKT0} \sup_{t\in(0,T)}\iidx{\Og}{|u|^2}+c_\ag\itQ{Q}{|u||Du|^2}\le C\itQ{Q}{|u|^2} + \iidx{\Og}{|u_0|^2}.\eeq Dropping the nonnegative second term on the left of \mref{keyestSKT0}, we obtain an integral Gronwall inequality for $y(t)=\|u(t)\|_{L^2(\Og)}^2$ in $(0,T_0)$ so that $\|u(t)\|_{L^2(\Og)}^2$ is bounded by a constant $c(T_0,\|u_0\|_{L^2(\Og)})$. This proves the lemma. \eproof

We now consider the general case and consider the condition \beqno{SKTcoop} \myprod{f(u),u} \le C_0|u|^2+c_0|u|^3,\;\forall u\in \RR^2.\eeq

\blemm{CufuSKTlem1} Assume \mref{SKTcoop}. Then the conclusion of \reflemm{CufuSKTlem} still holds if either that homogeneous Dirichlet boundary condition is assumed and $c_0$ is small or that homogeneous Neumann boundary condition is assumed and $\sup_{(0,T_0)}\|u\|_{L^1(\Og)}$ is bounded. \elemm

\bproof Revisiting the proof of \reflemm{CufuSKTlem}, we need only show that a similar version of \mref{keyestSKT0} holds here to give a Gronwall inequality for $\|u(t)\|_{L^2(\Og)}$ so that the proof of \reflemm{CufuSKTlem} can continue.
Indeed, instead of \mref{keyestSKT0} we now have for $T\in(t_0,T_0)$ \beqno{keyestSKT1} \sup_{t\in(t_0,T)}\iidx{\Og}{|u|^2}+c_\ag\itQ{Q}{|u||Du|^2}\le \itQ{Q}{(C_0|u|^2+c_0|u|^3)}+\iidx{\Og}{|u_0|^2}\eeq

If $u=0$ on $\partial\Og\times (0,T_0)$ then an application of Poincar\'e's inequality to $|u|^{3/2}$ yields
$$\iidx{\Og}{|u|^3} \le C\iidx{\Og}{|u||Du|^2}$$ for some constant $C$ depending only on $N$. Thus, if $c_0$ is small in terms of $c_\ag$ then the integral of $u^3$ in \mref{keyestSKT1} can be absorbed into the left hand side so that we obtain \mref{keyestSKT0}. 

Otherwise, for Neumann boundary condition, by \reflemm{Sobointineq} we see that for any given $\eg>0$ and $\bg\in(0,1]$ there is a constant $C(\eg,\bg)$ such that
$$\iidx{\Og}{|u|^3} \le \eg \iidx{\Og}{|u||Du|^2}+C(\eg,\bg)\left(\iidx{\Og}{|u|^{\frac32\bg}}\right)^{\frac2\bg}.$$

We now choose $\bg=2/3$ and $\eg$ sufficiently small to see that \mref{keyestSKT1} gives
\beqno{keyestSKT2} \sup_{t\in(t_0,T_0)}\iidx{\Og}{|u|^2}+\frac{c_\ag}{2}\itQ{Q}{|u||Du|^2}\le C\itQ{Q}{|u|^2}+C\sup_{(0,T_0)}\|u\|_{L^1(\Og)}^3.\eeq This is similar to \mref{keyestSKT0} and if $\sup_{(0,T_0)}\|u\|_{L^1(\Og)}$ is bounded then we obtain again a Gronwall inequality like \mref{keyestSKT0}. The proof is complete. \eproof

\section{Uniqueness of limiting solutions} \label{uniweaksec}\eqnoset
\newcommand{\mK}{\mathbf{K}}

We discuss the uniqueness of weak solutions obtained as limits of strong solutions in  the approximation process described in \refsec{weaksol}. We will show that any subsequence of these strong solutions in fact converges to a unique weak solution. As a consequence, the whole sequence converges to this limiting weak solution. Similar results for more general approximation schemes will be discussed in \refrem{BigPrem}.

\btheo{uniKthm} Suppose that $f(u)=Ku + g(u)$ for some constant $m\times m$ matrix $K$ and $g$ satisfies $|g_u(u)|\lesssim \llg (u)$ for all $u\in\RR^m$. Then the weak solution obtained by the approximation process in \reftheo{strongweakthm} is unique. \etheo

\bproof We consider two approximation schemes  with $P_{i,n}(u)=\llg_{i,n}u+P(u)$ ($i=1,2$) for some sequences $\llg_{i,n}\to0$
\beqno{approx12}\left\{\barr{ll}u_t=\Delta (P_{i,n}(u)) +f(u) &\mbox{in $\Og \times(0,T_0)$},\\ u=0 &\mbox{on $\partial\Og \times (0,T_0)$},\\u=u_0&\mbox{ on $\Og $}. \earr \right.\eeq

Let $\{u_{1,n}\}$ and $\{u_{2,n}\}$ be the sequences of strong solutions of \mref{approx12} that converge to the two weak solutions $u_1,u_2$ respectively. We will show that $u_1\equiv u_2$ on $\Og \times(0,T_0)$.

For any integers $m,n$, subtracting the equations of $u_{1,n}$ and $u_{2,m}$, we get for $w:=u_{1,n}-u_{2,m}$
\beqno{weqnstart}w_t=\Delta(P_{1,n}(u_{1,n})-P_{2,m}(u_{2,m}))+f(u_{1,n})-f(u_{2,m}).\eeq

We can write
$$f(u_{1,n})-f(u_{2,m})= Kw+\mG_{m,n} w,$$
where we denoted $$\mG_{m,n}=\int_0^1\frac{\partial}{\partial u}g(su_{1,n}+(1-s){u_{2,m}})\,ds.$$

Also,
$$\barr{lll}P_{1,n}(u_{1,n})-P_{2,m}(u_{2,m})&=&P_{1,n}(u_{1,n})-P_{1,n}(u_{2,m})+P_{1,n}(u_{2,m})-P_{2,m}(u_{2,m})\\&=&\mA_{m,n}^{(1)} w+P_{1,n}(u_{2,m})-P_{2,m}(u_{2,m}),\earr$$
$$\mA_{m,n}^{(1)}=\int_0^1\frac{\partial}{\partial u}P_{1,n}(su_{1,n}+(1-s){u_{2,m}})\,ds.$$
Similarly,
$$\barr{lll}P_{1,n}(u_{1,n})-P_{2,m}(u_{2,m})&=&P_{1,n}(u_{1,n})-P_{2,m}(u_{1,n})+P_{2,m}(u_{1,n})-P_{2,m}(u_{2,m})\\&=&P_{1,n}(u_{1,n})-P_{2,m}(u_{1,n})+\mA_{m,n}^{(2)} w,\earr$$
$$\mA_{m,n}^{(2)}=\int_0^1\frac{\partial}{\partial u}P_{2,m}(su_{1,n}+(1-s){u_{2,m}})\,ds.$$

Define $\mA_{m,n}=\frac12(\mA_{m,n}^{(1)}+\mA_{m,n}^{(2)})$ and $$\mP_{m,n}=\frac12(P_{1,n}(u_{2,m})-P_{2,m}(u_{2,m})+P_{1,n}(u_{1,n})-P_{2,m}(u_{1,n})).$$ Using these introduced terms in \mref{weqnstart} we easily see that
\beqno{weqnstart1}w_t=\Delta(\mA_{m,n}w)+\Delta\mP_{m,n}+K w+\mG_{m,n}w.\eeq

Hence,  for any $T\in (0,T_0)$ and  $\Psi$ on $L^2(\Og\times(0,T))$
$$\itQ{Q(s)}{\myprod{w_t,\Psi}}=\itQ{Q(s)}{\myprod{\Delta(\mA_{m,n} w)+(\mG_{m,n}+K) w,\Psi}}+\itQ{Q(s)}{\myprod{\Delta\mP_{m,n},\Psi}},$$ where we denoted  $Q(s)=\Og\times (0,s)$ for any  $s\in(0,T)$.

Assume that $\Psi$ is sufficiently smooth and satisfies $\Psi=0$ on $\partial\Og \times (0,T)$. Integrating by parts twice in $x$ ($\Psi, u_{1,n}, u_{2,m}$ are zero on the boundary)  and rearranging, we have
\beqno{wneqn}\barr{lll}\itQ{Q(s)}{\myprod{w,\Psi}_t}&=&\itQ{Q(s)}{\myprod{\Psi_t+\mA_{m,n}^T\Delta \Psi+\mG_{m,n}^T\Psi+K^T\Psi, w}}\\&&+\itQ{Q(s)}{\myprod{ \mP_{m,n}, \Delta\Psi}}.\earr\eeq

Concerning the first integral on the right hand side,  \reflemm{Psiexistence} following this proof shows that for any given $\psi\in C^{1}(\Og)$ there is a sequence of strong solutions $\Psi_{m,n}$ to
\beqno{Psidefa0}\left\{\barr{ll}\Psi_t+\mA_{m,n}^T\Delta \Psi+\mG_{m,n}^T\Psi+K^T\Psi=0&\mbox{on $Q=\Og \times(0,T)$},\\\Psi=0&\mbox{on $\partial\Og \times(0,T)$},\\ \Psi(x,T)=\psi(x).&\earr\right.\eeq 
We will also show in  \reflemm{Psiexistence} that there is a constant $C(T_0,\|\psi\|_{C^1(\Og)})$ such that 
\beqno{Psibound} \itQ{Q}{\llg_{(m,n)}|\Delta\Psi_{m,n}|^2} \le C(T_0,\|\psi\|_{C^1(\Og)}),\eeq where $\llg_{(m,n)}$ is the ellipticity function of the matrix $\mA_{m,n}$. 

Combining \mref{wneqn} and \mref{Psidefa0}, because $w(0)=0$, we have
\beqno{Psimn}\iidx{\Og\times\{s\}}{w\Psi_{m,n}}= I_{m,n}(s),\eeq
where, as $\mP_{m,n}=\frac12(\llg_{1,n}-\llg_{2,m})(u_{1,n}+u_{2,m})$, $$I_{m,n}(s)=\frac12\itQ{Q(s)}{(\llg_{1,n}-\llg_{2,m})\myprod{(u_{1,n}+u_{2,m}) ,\Delta\Psi_{m.n}}}.$$

Clearly, $|\llg_{1,n}-\llg_{2,m}|\le |\llg_{1,n}-\llg_{2,m}|^\frac12|\llg_{1,n}+\llg_{2,m}|^\frac12$. By H\"older's inequality we have 
$$\barrl{|I_{m,n}(s)|\le \frac12 |\llg_{1,n}-\llg_{2,m}|^\frac12\itQ{Q(s)}{ |u_{1,n}+u_{2,m}||\llg_{1,n}+\llg_{2,m}|^\frac12||\Delta\Psi_{m.n}|}}{.5cm}&\le \frac12|\llg_{1,n}-\llg_{2,m}|^\frac12\left(\itQ{Q}{|u_{1,n}+u_{2,m}|^2}\right)^\frac12\left(\itQ{Q}{|\llg_{1,n}+\llg_{2,m}||\Delta\Psi_{m,n}|^2}\right)^\frac12.\earr$$

By \mref{u2est0}, the first integral on the right hand side is bounded. On the other hand, it is clear that the ellipticity function $\llg_{(m,n)}$ of the matrix $\mA_{m,n}$  satisfies $\llg_{(m,n)}\ge \frac12(\llg_{1,n}+\llg_{2,m})$. Thus, by \mref{Psibound}, the second integral on the right hand side is also bounded. As $\llg_{1,n},\llg_{2,m}\to0$, we conclude that $I_{m,n}(s)\to0$ as $m,n\to\infty$.

Let $\psi_0\in C^1(\Og \times (0,T_0))$. For any $T\in(0,T_0)$ we take $\psi=\psi_0(T)$ in the above argument. We just showed that the functions  $$W_{m,n}(s):=\iidx{\Og\times\{s\}}{\myprod{u_{1,n}-u_{2,m},\psi}}=\iidx{\Og\times\{s\}}{\myprod{u_{1,n}-u_{2,m},\psi_{0}}}$$ converges to 0 on $(0,T_0)$. Using the fact that the $L^2(\Og)$ norms of $u_{1,n}(t),u_{2,m}(t)$ (see \mref{u2est}) are bounded uniformly on $(0,T_0)$, we see that $W_{m,n}$'s are also bounded uniformly on $(0,T_0)$. By the Dominated convergence theorem, we conclude that $W_{m,n}\to0$ in $L^1(0,T_0)$. Of course, as $u_{1,m}, u_{2,m}$ converge weakly to $u_1,u_2$ in $L^2(Q)$, we then have  $$ \itQ{Q}{\myprod{u_1-u_2,\psi_0}}=\lim_{(m,n)\to\infty}\int_0^{T_0}W_{(m,n)}(s)\,ds=0.$$

Because $C^1(\Og \times (0,T_0))$ is dense in $L^2(Q)$, the above also holds for all $\psi_0\in L^2(Q)$ and we see that $u_1\equiv u_2$ in $Q$.
\eproof

\brem{BigPrem} 
We can consider a more general approximation scheme by considering $P_{i,n}(u)=\pi_{i,n}(u)u+P(u)$ with $\pi_{i,n}(u)u$ being 'regularizers' in the sense that PR) holds for $P_{i,n}$. The same argument in \reftheo{strongweakthm} provides the existence of weak solutions. By the same proof of \reftheo{uniKthm}, we can prove the uniqueness result as long as we can establish its two key facts: The existence of the sequence $\Psi_{m,n}$ and that $I_{m.n}\to0$. The first one is easy because PR) is satisfied here. Concerning $I_{m,n}$, we replace $\llg_{i,n}$ in the proof by $\pi_{i,n}(u)$ to see that 
$$\mP_{m,n}=\frac12([\pi_{1,n}(u_{2,m})-\pi_{2,m}(u_{2,m})]u_{2,m}+[\pi_{1,n}(u_{1,n})-\pi_{2,m}(u_{1,n})]u_{1,n}).$$
Furthermore, we can also assume for $u=u_{i,n}$, $i=1,2$, that $|\pi_{1,n}(u)-\pi_{2,n}(u)|\le \llg_{(m,n)}$, the ellipticity constant of $\mA_{m,n}$, so that $I_{m,n}$ can be estimated by the integrals
$$\left(\itQ{Q}{|\pi_{1,n}(u_{i,n})-\pi_{2,m}(u_{i,n})||u_{i,n}|^2}\right)^\frac12\left(\itQ{Q}{|\llg_{(m,n)}||\Delta\Psi_{m,n}|^2}\right)^\frac12.$$ 

Hence, if $\|\pi_{1,n}(u)-\pi_{2,n}(u)\|_{L^\infty(\RR^m)} \to0$ for $u=u_{i,n}$ then the argument immediately goes through. Of course, if the norms $\|u_{i,n}\|_{L^{2q}(Q)}$'s are bounded for some $q\ge1$ then this condition can be improved by requiring only that $$\|\pi_{1,n}(u)-\pi_{2,n}(u)\|_{L^{q'}(\RR^m)} \to0 \mbox{ for $u=u_{i,n}$}.$$
\erem

We now present the key lemma providing the existence of the sequence $\{\Psi_{m,n}\}$ used in the proof of \reftheo{uniKthm}.
\blemm{Psiexistence} Let $\psi\in C^{1}(\Og,\RR^m)$. For any integers $m,n$, and $T\in(0,T_0)$  there is a function $\Psi_{m,n}$ solving
\beqno{Psidefa}\left\{\barr{ll}\Psi_t+\mA_{m,n}^T\Delta \Psi+\mG_{m,n}^T\Psi+K^T\Psi=0&\mbox{on $Q=\Og \times(0,T)$},\\\Psi=0&\mbox{on $\partial\Og \times(0,T)$},\\ \Psi(x,T)=\psi(x).&\earr\right.\eeq 

In addition, there is a constant $C(T,\|\psi\|_{C^1(\Og)})$  such that
\beqno{Psibounda} \iidx{\Og\times\{s\}}{|D\Psi_{m,n}|^2} \le C(T,\|\psi\|_{C^1(\Og)}) \mbox{  for all $s\in[0,T]$},\eeq
\beqno{D2Psibound} \itQ{Q}{\llg_{(m,n)}|\Delta\Psi_{m,n}|^2} \le C(T,\|\psi\|_{C^1(\Og)}),\eeq where $\llg_{(m,n)}$ is the ellipticity function of the matrix $\mA_{m,n}$.
\elemm

\bproof Using a change of variables $t\to T-t$ the system \mref{Psidefa} is equivalent to the following linear parabolic system with homogeneous Dirichlet boundary condition and initial data $\hat{\Psi}(x,0)=\psi(x)$ for $\hat{\Psi}(x,t)=\Psi(x,T-t)$.
\beqno{Psihatdef}\hat{\Psi}_t=\mA_{m,n}^T(x,T-t)\Delta \hat{\Psi}+\mG_{m,n}^T(x,T-t)\hat{\Psi}+K^T\hat{\Psi}\mbox{ on $\Og \times(0,T)$},\eeq

Because $u_{1,n}$ and $u_{2,m}$ are strong solutions the coefficients of the above systems are smooth and bounded. From P) and the definition of $\mA_{m,n}$ we see that the ellipticity function $\llg_{m,n}$ of $\mA_{m,n}$  satisifies $\llg_{m,n}\ge\frac12(\llg_{1,n}+\llg_{2,m})>0$. Hence, the above system is a regular linear parabolic system so that it has a strong solution $\hat{\Psi}$. Thus, $\Psi_{m,n}$ exists.

We temporarily drop the subscripts $m,n$ in the calculation below. Multiplying \mref{Psidefa} with $\Delta\Psi$ and integrating by parts ($\Psi_t=0$ on the boundary because $\Psi$ is), we get for any $s<T'<T$ and $Q^{(s)}=\Og\times(s,T')$ \beqno{Psistart}-\itQ{Q^{(s)}}{\frac{d}{dt}|D\Psi|^2}+\itQ{Q^{(s)}}{\myprod{\mA^T\Delta\Psi,\Delta\Psi}}=-\itQ{Q^{(s)}}{\myprod{(\mG^T+K^T)\Psi,\Delta\Psi}}.\eeq

By P), for any vector $\zeta$ we can find a positive function $\llg_*$ such that
$$\myprod{\mA(u,v)\zeta,\zeta}\ge \int_0^1\llg(su+(1-s)v)\,ds|\zeta|^2\Rightarrow \myprod{\mA^T\Delta\Psi,\Delta\Psi}\ge \llg_*|\Delta\Psi|^2.$$

We now estimate the integral on the right hand side of \mref{Psistart}. First of all, integrating by parts in $x$, we have
$$-\itQ{Q^{(s)}}{\myprod{K^T\Psi,\Delta\Psi}}=\itQ{Q^{(s)}}{\myprod{K^TD\Psi,D\Psi}}\le C\itQ{Q^{(s)}}{|D\Psi|^2}.$$ 

Next, from the growth assumption $|g_u(u)|\le C\llg (u)$ and the definition of $\mG $, we see that $|\mG|\le C\llg_*$. So that by Young's inequality  
$$\myprod{\mG^T\Psi,\Delta\Psi}\le \eg\llg_*|\Delta\Psi|^2+C(\eg)\llg_*|\Psi|^2.$$ 

Therefore, for small $\eg>0$ we deduce from the above estimates and \mref{Psistart} the following inequality 
\beqno{Psistartz}\barrl{\iidx{\Og\times\{s\}}{|D\Psi|^2}+\itQ{Q^{(s)}}{\llg_*|\Delta\Psi|^2}\le}{2cm} &\iidx{\Og\times\{T'\}}{|D\Psi|^2}+C\itQ{Q^{(s)}}{\llg_*|\Psi|^2}+C\itQ{Q^{(s)}}{|D\Psi|^2}.\earr\eeq

Choosing $q$ such that $N/2<q<q_0$, we easily see that $2q'<2_*$ so that we can  estimate the integral of $\llg_{*}|\Psi|^2$ over $\Og $ by, using  H\"older and Sobolev's inequalities ($\Psi=0$ on the boundary) \beqno{DPsibound} \left(\iidx{\Og\times\{\tau\}}{\llg_*^{q}}\right)^\frac1{q}\left(\iidx{\Og\times\{\tau\}}{|\Psi|^{2q'}}\right)^\frac1{q'}\le C\iidx{\Og\times\{\tau\}}{|D\Psi|^{2}},\quad \tau\in(0,T).\eeq
Here, we used the fact that $\llg_{*}$ satisfies the same uniform bound of the ellipticity function of the matrix $P_{u}(u_{i,n})$ in \reftheo{strongweakthm}. Namely, $\|\llg(u_{i,n})\|_{L^{q_0}(\Og)}\le C$ for $i=1,2$ and some constant $C$, so that  \beqno{llg*bound}\|\llg_*\|_{L^{q_0}(\Og\times\{\tau\})}\le C \quad \forall \tau\in (0,T).\eeq Hence, $$\itQ{Q^{(s)}}{\llg_*|\Psi|^2}\le C\itQ{Q^{(s)}}{|D\Psi|^2}.$$ Using this in \mref{Psistartz} we deduce \beqno{Psistart0}\iidx{\Og\times\{s\}}{|D\Psi|^2}+\itQ{Q^{(s)}}{\llg_*|\Delta\Psi|^2}\le \iidx{\Og\times\{T'\}}{|D\Psi|^2}+C\itQ{Q^{(s)}}{|D\Psi|^2}.\eeq

By \refrem{limDuat0} after this proof, we have  $$\liminf_{t\to T}\|D\Psi(\cdot,t)\|_{L^2(\Og)}=\liminf_{t\to0}\|D\hat{\Psi}(\cdot,t)\|_{L^2(\Og)}\le C(\|\psi\|_{C^1(\Og)}).$$ Let $\{T_k\}$ be a sequence such that $T_k<T$ and $$\lim_{k\to\infty}\|D\Psi(\cdot,T_k)\|_{L^2(\Og)}=\liminf_{t\to T}\|D\Psi(\cdot,t)\|_{L^2(\Og)}.$$ Replacing $T'$ in \mref{Psistart0} by $T_k$ and letting $k\to \infty$, we then obtain
\beqno{Psistart1}\iidx{\Og\times\{t\}}{|D\Psi|^2}+\itQ{Q^{(s)}}{\llg_*|\Delta\Psi|^2}\le C(\|\psi\|_{C^1(\Og)})+C\int_s^T\iidx{\Og\times\{s\}}{|D\Psi|^2}\,dt.\eeq

This is an integral  Gronwall inequality for $\|D\Psi\|_{\Og\times\{s\}}^2$ which yields $\|D\Psi\|_{\Og\times\{s\}}\le C(T,\|\psi\|_{C^1(\Og)})$ for some constant
$C(T,\|\psi\|_{C^1(\Og)})$.
This is \mref{Psibounda}.
 
We also obtain the estimate \mref{D2Psibound} for $\Delta\Psi$ from \mref{Psistart1} and \mref{Psibounda}.
This completes the proof of the lemma. \eproof

\brem{limDuat0} In the proof, we used a result that $\liminf_{t\to0}\|D\hat{\Psi}\|_{L^2(\Og \times\{t\})}$ is bounded by some constant depending on $\|\psi\|_{C^1(\Og)}$ for solution of \mref{Psihatdef} with initial data $\psi$. This fact is in the same spirit of the Hille-Yoshida theorem (e.g., see \cite[Theorem 7.8]{brezis} or \cite[Theorem 5]{Evans}) concerning the continuity of $\|D\hat{\Psi}\|_{L^2(\Og \times\{t\})}$ when $\mA_{m,n}$ is a constant or independent of $t$. The matter is a bit subtle otherwise. More importantly, the estimate for $\liminf_{t\to0}\|D\hat{\Psi}\|_{L^2(\Og \times\{t\})}$ should {\em not} depend on higher order norms of $\mA_{m,n},\mG_{m,n}$ although that they are smooth. Also, this estimate must be uniform or independent of the ellipticity constant $\llg_{(m,n)}$ as it will tend to 0.  
As we cannot find an appropriate reference for this fact, we sketch the proof here. 

We split $\hat{\Psi}=h+H$ where $h,H$ solve 
$$h_t=a(x) \Delta h +\ccG h,\; H_t= \ccA \Delta H +\ccB \Delta h+ \ccG H,$$ where $\ccA(x,t)=\mA_{m,n}^T(x,T-t)$ , $\ccG(x,t)=\mG_{m,n}^T(x,T-t)$, $a(x)=\ccA(x,0)-\frac12\llg_{*,0}$ ($\llg_{*,0}=\llg_{(m,n)}$, the ellipticity constant for $\mA_{m,n}^T$) and $\ccB(x,t)=\ccA(x,t)-a(x)$. Also, $h(0)=\psi$ and $H(0)=0$. We rewrite the equation for $H$ as
$$H_t= \Div(\ccA DH +\ccB Dh)- D\ccA DH-D\ccB Dh + \ccG H$$
and test the system with $H$ and use the fact that $H(0)=0$ to obtain  for any $s>0$ that $$\barrl{\iidx{\Og \times\{s\}}{|H|^2}+\int_0^s\iidx{\Og}{\myprod{\ccA DH,DH}}=}{3cm}&-\dspl{\int_0^s}\iidx{\Og}{(\myprod{\ccB Dh,DH}+\myprod{D\ccA DH+D\ccB Dh,H}+\myprod{\ccG H,H})}.\earr$$  Applying Young's inequalities to the integrals on the right hand side  and using the ellipticity of $\ccA$, we easily get
$$\barrl{\llg_{*,0}\int_0^s\iidx{\Og}{|DH|^2}dt \le C\int_0^s\iidx{\Og}{|\ccB|^2\llg_{*,0}^{-1} |Dh|^{2}}dt+}{3cm}& C\dspl{\int_0^s}\iidx{\Og}{(|D\ccA|^2+|\ccG|)|H|^2}dt+C\dspl{\int_0^s}\iidx{\Og}{|D\ccB|Dh|H|}dt.\earr$$
We now divide the about inequality by $s$ to have
\beqno{keyDuat0}\barrl{\llg_{*,0}\frac1s\int_0^s\iidx{\Og}{|DH|^2}dt \le C\frac1s\int_0^s\iidx{\Og}{|\ccB|^2\llg_{*,0}^{-1} |Dh|^{2}}dt+}{2cm}& C \dspl{\frac1s\int_0^s}\iidx{\Og}{(|D\ccA|^2+|\ccG|) |H|^2}dt+C\frac1s\dspl{\int_0^s}\iidx{\Og}{|D\ccB|Dh|H|}dt.\earr\eeq
We will let $s\to0$ and need to investigate the limits of the terms on the right hand side.

From the definition of $\ccB$, $\ccB(x,t)=\ccA(x,t)-\ccA(x,0)+\frac12\llg_{*,0}$. By the continuity of $\ccA $ at $0$, we see that $\lim_{t\to0}|\ccB|^2\llg_{*,0}^{-1}\sim \llg_{*,0}$. By the Hille-Yoshida theorem (e.g., see \cite[Theorem 7.8]{brezis}), note that $a(x)$ is elliptic, smooth and independent of $t$, $Dh$ belongs to $C([0,T_0],L^{2}(\Og))$. In particular, $\|Dh(t)\|_{L^{2}(\Og)}$ is  continuous at $t=0$. Hence, the limit of the first term on the right hand side of \mref{keyDuat0} when $s\to0$ is bounded by a multiple of  $\llg_{*,0}\|\psi\|_{C^1(\Og)}$.

Meanwhile,  $D\ccA$ and $D\ccB$ are bounded near $t=0$  because they depend on the spatial derivatives of the strong solutions $u_{1,n},u_{2,m}$ at $t=T$. As $H(0)=0$, the last two terms tend to $0$. 

Hence, letting $s\to0$ in \mref{keyDuat0}, we derive $\llg_{*,0}\liminf_{t\to0}\|DH\|_{L^2(\Og \times\{t\})}\le C\llg_{*,0}\|\psi\|_{C^1(\Og)}$ so that
$\liminf_{t\to0}\|DH\|_{L^2(\Og \times\{t\})}\le C(\|\psi\|_{C^1(\Og)})$. As $\hat{\Psi}=h+H$, we obtain the desired bound for $\liminf_{t\to0}\|D\hat{\Psi}\|_{L^2(\Og \times\{t\})}$. 
\erem

\section{Uniqueness of (unbounded) weak solutions} \eqnoset\label{uniboundw}
We have proved that the weak solution obtained by the limiting process in the \refsec{weaksol} is unique. To the best of our knowledge, the existence of weak solutions to the degenerate scalar equations has always been established by this way in literature. It is desirable to establish the uniqueness of general weak solutions, whose existence can be established by different methods. This has been done for {\em bounded} weak solutions of {\em scalar} equations/systems (even if they are degenerate in some cases \cite{Vas}). But this is not a satisfactory result for systems because the boundedness of solutions to systems generally is an open problem and the arguments for scalar equations are not applicable here. However, if the systems are nondegenerate then we can establish a uniqueness result for unbounded weak solutions based on a similar argument in the proof of \reftheo{uniKthm}.

Following the definition W) in \refsec{weaksol} for weak solutions to the degenerate case, we say that $u$ is a weak solution on $\Og \times(0,T_0)$ of \mref{ep1} if for  a.e. $T\in(0,T_0)$ and any $\fg\in C^1(\Og \times (0,T))$ we have, provided that the following integrals are all finite
\beqno{wdef}\barrl{\iidx{\Og}{\myprod{u(T),\fg(T)}-\myprod{u_0,\fg(0)}}=}{3cm}&\itQ{\Og \times(0,T)}{[\myprod{u, \fg_t}-\myprod{A(u)Du,D\fg} +\myprod{f(u),\fg}]}.\earr\eeq

Clearly, in oder for the above integrals are finite for all $\fg\in C^1(\Og \times (0,T))$, we need to impose, at least, that $u\in L^\infty((0,T_0),L^1(\Og)$ and $A(u)Du\in L^1(\Og \times (0,T_0))$.

As  $A(u)=P_u(u)$ so $A(u)Du=D(P(u))$, the above equation easily gives
\beqno{Pwdef}\barrl{\iidx{\Og}{\myprod{u(T),\fg(T)}-\myprod{u_0,\fg(0)}}=}{3cm}&\itQ{\Og \times(0,T)}{[\myprod{u, \fg_t}+\myprod{P(u),\Delta\fg} + \myprod{f(u),\fg}]}.\earr\eeq

Our main result in this section states that if such a weak solution satisfies some very mild integrability conditions then it is unique.

\btheo{uniweak} Assume that $P_u$ is regular elliptic. That is there are function $\llg$ and constant $\llg_0>0$ such that $\llg (u)\ge\llg_0$ and $$ \myprod{P_u(u)\zeta,\zeta}\ge \llg(u) |\zeta|^2, \mbox{ for all $u\in\RR^m, \zeta\in\RR^{Nm}$}.$$

Assume that for some $p>2$ the maps $u\to \partial_u P(u)$ and $u\to \partial_u f(u)$ are continuous from $L^p(Q)$ to $L^q(Q)$, $Q=\Og \times(0,T_0)$ and $q=2p/(p-2)$.

If $u$ is a weak solution satisfying $u\in L^p(Q)$ and \beqno{llguw} \sup_{t\in (0,T_0)}\|\llg(u(t))\|_{L^{q_0}(\Og)}<\infty \mbox{ for some $q_0>N/2$},\eeq then $u$ is unique. \etheo

\newcommand{\ma}{\mathbf{a}}
\newcommand{\mg}{\mathbf{g}}

\bproof For any $u_1,u_2$ we can write  $$P(u_1)-P(u_2)=\ma(u_1,u_2)(u_1-u_2),\quad \ma(u_1,u_2):=\int_0^1 \partial_u P(su_1+(1-s)u_2)\,ds,$$ 
$$f(u_1)-f(u_2)=\mg(u_1,u_2)(u_1-u_2),\quad \mg(u_1,u_2):=\int_0^1 \partial_u f(su_1+(1-s)u_2)\,ds.$$

Using these notations, if $u_1,u_2$ are two weak solutions with the same initial data $u_0$ then we subtract the two systems \mref{Pwdef} for $u_1,u_2$ to see that $w=u_1-u_2$ satisfies
\beqno{u1u2bded}\iidx{\Og}{\myprod{w(T),\fg(T)}}=\itQ{\Og \times(0,T)}{\myprod{w, \fg_t+\ma(u_1,u_2)^T\Delta\fg + \mg(u_1,u_2)^T\fg}}.\eeq

We consider the sequences $\{u_{1,n}\}$, $\{u_{2,n}\}$ of mollifications of $u_1, u_2$. That is, we consider $C^\infty$ functions $\eta(t)$ and $\rg (x)$ whose supports are $(-1,1)$ and $B_1(0)$ and $\|\eta\|_{L^1(\RR)} =\|\rg\|_{L^1(\RR^N)}=1$. Denote $\eta_n(t) = n\eta(t/n)$ and $\rg_n(x)=n^N\rg (x/n)$.
For $i=1,2$ define $$ u_{i,n}(t,y)=(\eta_n\fg_n)* u_i(t,y)=\int_{\RR}\iidx{\RR^N}{\eta_n (s-t)\fg_n(x-y)u_i(t,x)}\,ds.$$

We now follow a similar argument as in the proof of \reftheo{uniKthm}. For each $n$ the matrix $\ma(u_{1,n},u_{2,n})$ is smooth and uniformly elliptic with the elliptic function
$$\hat{\llg}_n(x,t) = \int_0^1 \llg(su_{1,n}+(1-s)u_{2,n})\,ds.$$
Hence,  for any $\psi\in C^{1}(\Og)$ and $T\in(0,T_0)$ \reflemm{Psiexistence}, with the growth assumptions on $g_d(u_{1,n},u_{2,n})$,  provides  strong  (classical) solutions $\Psi_{n}$ to the systems
\beqno{Psisystem}\left\{\barr{l}\Psi_t+\ma(u_{1,n},u_{2,n})^T\Delta\Psi + \mg(u_{1,n},u_{2,n})^T\Psi=0\mbox{ in $Q:=\Og \times (0,T)$},\\\Psi=0 \mbox{ on $\partial\Og \times (0,T)$},\\\Psi(x,T) = \psi(x) \mbox{ on $\Og $}.\earr \right.\eeq

Furthermore, because $\llg $ is convex, by Jensen's inequality $\hat{\llg}_n\le \frac12(\llg (u_{1,n})+\llg (u_{2,n}))$. Similarly, for $i=1,2$ $\llg (u_{i,n})\le (\eta_n\rg_n) *\llg (u_i)$ so that
$$\|\llg (u_{i,n}(t))\|_{L^{q_0}(\Og)} \le \int_\RR \eta_n(s-t)\|\rg_n *_x \llg(u_{i}(t))\|_{L^{q_0}(\Og)}\,ds.$$
Here, $*_x$ denotes the convolution in $\RR^N$. Because $\|\eta_n\|_{L^1(\RR)}=1$ and $\|\rg_n *_x \llg(u_{i}(t))\|_{L^{q_0}(\Og)}\le\|\llg (u_{i}(t))\|_{L^{q_0}(\Og)} \le C$ by the assumption \mref{llguw}, we have 
\beqno{llgun} \|\llg (u_{i,n}(t))\|_{L^{q_0}(\Og)} \le C \mbox{ for all $t\in(0,T_0)$ and integer $n$}.\eeq

Hence, the ellipticity functions $\hat{\llg}_n$ of $\ma(u_{1,n},u_{2,n})$ satisfies $\|\hat{\llg}_n\|_{L^{q_0}(\Og\times\{t\})}\le C$  on $(0,T_0)$ and we can apply  \reflemm{Psiexistence} to obtain a constant $C(\|\psi\|_{C^1(\Og)})$ such that for all $n$ \beqno{D2Wn}  \sup_{(0,T_0)}\|\Psi_n\|_{L^2(\Og)},\; \|\hat{\llg}_n^\frac12\Delta\Psi_n\|_{L^2(Q)}\le C(\|\psi\|_{C^1(\Og)}).\eeq

Let $\fg =\Psi_n$ in \mref{u1u2bded}, this is allowable  because $\Psi_n$ is a classical solution. From the equation of $\Psi_n$, we obtain
$$\barr{lll}\iidx{\Og}{\myprod{w(T),\psi}}&=&-\itQ{Q}{\myprod{w,[\ma(u_{1,n},u_{2,n})^T-\ma(u_{1},u_{2})^T]\Delta\Psi_n}}\\&&-\itQ{Q}{\myprod{w,[\mg(u_{1,n},u_{2,n})^T-\mg(u_{1},u_{2})^T]\Psi_n}}.\earr$$ This is \beqno{keyuniw}\barr{lll}\iidx{\Og}{\myprod{w(T),\psi}}&=&-\itQ{Q}{\myprod{[\ma(u_{1,n},u_{2,n})-\ma(u_{1},u_{2})]w,\Delta\Psi_n}}\\&&-\itQ{Q}{\myprod{[\mg(u_{1,n},u_{2,n})-\mg(u_{1},u_{2})]w,\Psi_n}}.\earr\eeq

Letting $n\to\infty$, we will see that the integrals on the right hand side tend to zero. Indeed, we consider the first integral. Because $\hat{\llg}_n\ge\llg_0>0$, the bound \mref{D2Wn} implies $\|\Delta\Psi_n\|_{L^2(Q)}$ is bounded uniformly so that  we need only to show that $[\ma(u_{1,n},u_{2,n})-\ma(u_{1},u_{2})]w$ converges strongly to 0 in $L^2(Q)$. By H\"older's inequality with $q=2p/(p-2)$
$$\|[\ma(u_{1,n},u_{2,n})-\ma(u_{1},u_{2})]w\|_{L^2(Q)}\le \|\ma(u_{1,n},u_{2,n})-\ma(u_{1},u_{2})\|_{L^q(Q)}\|w\|_{L^p(Q)}.$$

As we are assuming that the map $u\to \partial_u P(u)$ is continuous from $L^p(Q)$ to $L^q(Q)$ and because $u_{i,n}\to u_i$ in
$L^p(Q)$, it is clear from the definition of $\ma$ that $\ma(u_{1,n},u_{2,n})\to \ma(u_{1},u_{2})$ in $L^q(Q)$. Thus, $[\ma(u_{1,n},u_{2,n})-\ma(u_{1},u_{2})]w$ converges strongly to 0 in $L^2(Q)$. Thus, the first integral on the right hand side of \mref{keyuniw} tends to 0 as $n\to \infty$. 

Similar argument applies to the second integral to obtain the same conclusion. 
We just prove that the right hand side of \mref{keyuniw} tends to 0.
We then have $$\iidx{\Og}{\myprod{w(T),\psi}}=0 \mbox{ for any $\psi\in C^{1}(\Og)$}.$$ We conclude that $w(T)=0$ for all $T\in(0,T_0)$. Hence $u_1\equiv u_2$ on $Q$. \eproof

\brem{boundedwrem} If we discuss {\em bounded} weak solutions then the proof is much simpler as the needed convergences are obvious. As the sequences $\{u_{1,n}\}$, $\{u_{2,n}\}$  converge to $u_1,u_2$ in $L^\infty(Q)$ we see that $\ma(u_{1,n},u_{2,n})\to\ma(u_{1},u_{2}) $ and $\mg(u_{1,n},u_{2,n})\to\mg(u_{1},u_{2}) $ strongly in $L^\infty(Q)$. Even in this case, the result stated here is nontrivial because our class of weak solutions is much more larger than those considered in literature because our  class of  admissible test functions in \mref{wdef} is more restrictive, provided that it includes classical solutions of \mref{Psisystem}. Otherwise, as done in several works, one is allowed to take $\fg =u_i$ (a weak solution), $i=1,2$ in \mref{wdef} and subtract the results to obtain a Gronwall's inequality for $y(t):=\|u_1-u_2\|_{L^2(\Og \times\{t\})}$ and easily conclude that $y(t)\equiv0$ because $y(0)=0$ so that $u_1\equiv u_2$. \erem

The following consequence of \reftheo{uniweak} applying to cross diffusion models in mathematical biology models with polynomial growth data, the (SKT) system is an example.

\bcoro{ALqcont} Assume that $\partial_u P(u)$ and $\partial_u f(u)$ have polynomial growths
$$|\partial_u P(u)|,\; |\partial_u f(u)|\le C(|u|^k+1) \mbox{ for some $k>0$}.$$ Then the uniqueness conclusion of \reftheo{uniweak} applies to weak solutions in the space $L^p(Q)\cap L^\infty((0,T_0),L^r(\Og))$ if $p\ge 2(1+k)$ and $r>kN/2$.
\ecoro

\bproof We need only verify the assumptions of \reftheo{uniweak}. It is clear from its proof that we need to establish the convergences $\ma(u_{1,n},u_{2,n})\to \ma(u_{1},u_{2})$ and $\mg(u_{1,n},u_{2,n})\to \mg(u_{1},u_{2})$ in $L^q(Q)$ for $q=2p/(p-2)$ along some subsequences of $\{u_{1,n}\}$,$\{u_{2,n}\}$ which converge to $u_1$,$u_2$ in $L^p(Q)$ and a.e. in $Q$. By the Riesz-Fisher theorem we can find subsequences of $\{u_{i,n}\}$ and functions $\hat{u}_i\in L^p(Q)$ such that, after relabeling, $u_{i,n}\to u_i$ and $|u_{i,n}|\le \hat{u}_i$  a.e. in $Q$. The growth condition of $\partial_u P(u)$  then implies $|\ma(u_{1,n},u_{2,n})|\le |\hat{u}_1|^k+|\hat{u}_2|^k+1$, a function in $L^{p/k}(Q)$. Furthermore, $\ma(u_{1,n},u_{2,n})\to \ma(u_{1},u_{2})$ a.e. in $Q$ because $\ma$ is continuous. By the Dominated convergence theorem, we see that $\ma(u_{1,n},u_{2,n})\to \ma(u_{1},u_{2})$ in $L^{p/k}(Q)$. This also yields the convergence in $L^q(Q)$ for $q=2p/(p-2)$ because $q\le p/k$ as $p\ge 2+2k$. Similarly, from the growth condition of $\partial_u f$, we have
$\mg(u_{1,n},u_{2,n})\to \mg(u_{1},u_{2})$ in $L^q(Q)$.

Also, from the growth assumption on $\partial_u P$, we have that $|\llg(u)|^{q_0}\le C|u|^{kq_0}$. From the assumption of the corollary, $u\in L^\infty((0,T_0),L^r(\Og))$ for some $r>kN/2$ so that the quantity $\sup_{(0,T_0)}\|\llg (u)\|_{L^{q_0}(\Og)}$ is finite for some $q_0>N/2$.  We see that all assumptions of of \reftheo{uniweak} are verified here. This completes the proof.
\eproof

The space of test functions $\fg $ in our definition \mref{wdef} is the smallest possible one, $\fg\in C^1(\Og \times (0,T_0))$,  so that our class of weak solutions considered here is very wide, just sufficient for its integrals to be finite.
If we consider the definition of {\em generalized} solutions in \cite{LSU}, which has been commonly adopted to the definition weak solutions in many works, then we can apply the above results to many models in mathematical biology, including the SKT system. 

Following \cite[Chapter III]{LSU}, we say that $u$ is a generalized solution from $V_2(Q)$, the  Banach space with norm
$$\|u\|_{V_2(Q)}=\sup_{(0,T_0)}\|u\|_{L^2(\Og\times \{t\})}+\|Du\|_{L^2(Q)},$$
if $u$ satisfies \mref{wdef} for any $\fg \in W_2^{1,1}(Q)$, the Hilbert space with scalar product $$ \myprod{u,v}_{W_2^{1,1}(Q)}=\itQ{Q}{[\myprod{u,v}+\myprod{u_t,v_t}+\myprod{Du,Dv}]}.$$

 We also have the following inequality which can be proved easily by using Sobolev's embedding inequality in the same way as in \cite[(3.2) p.74]{LSU}.
For  any time interval $I$ and any nonegative measurable functions $g\in L^\infty(I, L^2(\Og))$, $G\in L^2(I, W^{1,2}(\Og))$  there is a constant $C$ such that\beqno{paraSobo}\itQ{\Og\times I}{g^{2r}G^2}\le C\sup_I\left(\iidx{\Og\times\{t\}}{g}\right)^{r}\itQ{\Og\times I}{(|DG|^2+G^2)},\eeq where $r=2/N$ if $N>2$ and $r$ is any number in $(0,1)$ if $N\le 2$.

Now, it is clear that in order for the integrals in \mref{wdef} are finite for all $\fg \in W_2^{1,1}(Q)$ we must assume further that $u\in L^\infty((0,T),L^2(\Og))$ and $D(P(u))\in L^2(Q)$. We now let $g=|u|$ and $G=|P(u)|$ in \mref{paraSobo}. Because $u\in V^2(Q)$ and $|P(u)|\sim |u|^{k+1}$, we see that $u\in L^{2r+2k+2}(Q)$. The condition in \refcoro{ALqcont} that $u\in L^{p}(Q)$ with $p\ge 2+2k$ is then obvious so that we need only that  $\sup_{(0,T_0)}\|u\|_{L^r(\Og)}$ is finite for some $r>kN/2$. This condition is clearly satisfied for generalized solutions ($r=2$) for the usual SKT system, where $P(u)$ has quadratic growth in $u$ (so that $k=1$), in domains with dimension $N\le3$. Our uniqueness result in \refcoro{ALqcont} then applies to this case and we can assert that if $u$ is a weak solution with $u\in L^\infty((0,T),L^2(\Og))$ and $D(P(u))\in L^2(Q)$ then $u$ is unique.

\bibliographystyle{plain}

\end{document}